\documentclass{article}

\usepackage[letterpaper,top=2cm,bottom=2cm,left=3cm,right=3cm,marginparwidth=1.75cm]{geometry}

\usepackage{amsmath,amssymb}
\usepackage{graphicx}
\usepackage[colorlinks=true, allcolors=blue]{hyperref}

\usepackage{amssymb, amsmath, amsthm, amsfonts, amscd}
\usepackage[english]{babel}
\usepackage{graphicx}
\usepackage{caption}
\usepackage{subcaption}
\usepackage{color,soul}
\usepackage{booktabs,dcolumn}                     
\usepackage[figureposition=bottom]{caption}       
\usepackage[usenames,dvipsnames,svgnames,table]{xcolor}
\usepackage{authblk}
\usepackage[title]{appendix}
\usepackage{tikz}
\usepackage{mathalfa}

\usepackage[utf8]{inputenc}
\usepackage[T1]{fontenc}
\usepackage{amsmath}
\usepackage{amsfonts}
\usepackage{amssymb}
\usepackage[version=4]{mhchem}
\usepackage{stmaryrd}
\usepackage{graphicx}
\usepackage[export]{adjustbox}

\usepackage{amsmath,amssymb,amsthm}
\usepackage{hyperref}
\usepackage{graphicx}
\usepackage{listings}
\usepackage{xcolor}
\usepackage[utf8]{inputenc}
\usepackage[T1]{fontenc}
\usepackage{amsmath,amsthm,amssymb,amsfonts}
\usepackage{geometry}

\newtheorem{thm}{Theorem}[section]
\newtheorem{lem}[thm]{Lemma}
\newtheorem{prop}[thm]{Proposition}
\newtheorem{cl}{Conjecture}[section]

\title{Extension of Boundary Control method to elliptic and parabolic problems, and its application to the Calderon problem}
\author{Dimitra Kyriakopoulou$^{\MakeLowercase{a}}$}
\affil{$^a$Mathematics Research Center, Academy of Athens, Greece}

\date{\vspace*{8mm}\textit{Dedicated to the living presence of Professor Yaroslav Kurylev}\vspace*{8mm}}

\begin{document}

\maketitle

\begin{abstract}
We show that Boundary Control method, a method for hyperbolic inverse problems, is also capable of dealing directly with certain classes of elliptic and parabolic Inverse Boundary Value Problems; thus pointing towards Boundary Control method potentially constituting a means of unification of Inverse Boundary Value Problems. 
As an application we show that the Calderon problem
can be dealt with directly via Boundary Control method, i.e. without reduction of the elliptic problem to a 'hyperbolized' problem. 
\end{abstract}

\section{Introduction} \label{Section_Intro}

Manifolds with
smooth boundary form a subcategory of manifolds with edge \cite{LSch}. 
Hence, Bernard Wolfgang Schulze studied elliptic Boundary Value Problems (BVPs) of Shapiro–
Lopatinskii (SL) type (\cite{SchInd},\cite{SchBPS},\cite{sing},\cite{mix}) in the framework
of Boutet de Monvel’s algebra of BVPs with the transmission property at the
boundary \cite{[10]}. Ellipticity refers to an analogue of Shapiro–
Lopatinskii conditions, i.e., a bijectivity condition for an operator-valued principal
symbol structure which contains also trace and potential operators with respect to
the edge, a substitute of the former boundary. 
Later Schulze \cite{LSch} created the Toeplitz analogue algebra of BVPs
unifying 
conditions of SL- and global projection (GP)- elliptic type (especially APS-conditions in the sense of Atiyah, Patodi, and Singer -\cite{[4]}, \cite{[5]}, \cite{[6]}); ellipticity in this context is equivalent
with the Fredholm properties in the respective scales of spaces (standard Sobolev
spaces in the SL case, spaces of Hardy type in the GP case). 

Reducing inverse elliptic and parabolic BVPs to edge problems, we will study how Boundary Control (BC) method can be applied to them. 
BC method, originating from M. G. Krein's work on 1-dimensional inverse scattering theory (\cite{[Kr51a]}, \cite{[Kr51b]}), contrasts with Gel’fand-Levitan and Marchenko's fundamental methods through its utilization of the wave equation's finite propagation speed. Krein's approach, while initially not evident due to its frequency domain formulation, was later clarified by Blagovestchenskii's time-domain analysis, highlighting its hyperbolic nature through a Volterra-type equation for unknown functions (\cite{[Bla71a]}). This advancement was essential for the method's multi-dimensional extension by Belishev \cite{[Be87]}, with its geometric aspects further detailed by Belishev and Kurylev \cite{belku}, also noted in \cite{[KKL01]}. Additionally, the BC method's capability to assess the inner product of boundary-induced waves, traced back to Blagovestchenskii (\cite{[Bla71b]}), was extended to multidimensional contexts (\cite{[BeBla92]}), underscoring the method's non-perturbative and inherently hyperbolic characteristics \cite{IsKu}.

The BC-method integrates control and systems theory, asymptotic methods, functional analysis, and operator theory with partial differential equations, indicating its interdisciplinary strength. It has established significant links with Banach algebras (\cite{bel17_[22]}, \cite{bel17_[62]}), non-commutative geometry (\cite{bel17_[41]}, \cite{bel17_[42]}, \cite{bel17_[37]}), and functional models of linear operators (\cite{bel17_[34]}, \cite{bel17_[56]}), underscoring the method's comprehensive applicability and depth \cite{bel17}.

Linking the BC method with noncommutative algebra enables its application to elliptic and parabolic problems, suggesting its potential as a unifying solution for Inverse Boundary Value Problems (BVPs). Noncommutative analysis, focused on quantizing observables within arbitrary Poisson brackets, aims to identify operators \(A_{1}, \ldots, A_{m}\) that express the problem's operator through a rich system of relations. This approach allows solving by expressing solutions as functions of these operators, where noncommutative analysis outlines conditions for this method and provides formulas for symbol composition laws. This reduction simplifies the approach to asymptotic problems, making noncommutative analysis particularly effective for various asymptotics. It outperforms wave-packet transform methods in handling simultaneous asymptotics, proving useful in constructing parametrices on manifolds with singularities, applying to asymptotics in weighted function spaces, such as power-law weighted Sobolev spaces near degeneration points \cite{Qm}.

We extend the BC method to noncommutative algebras, paralleling the commutative \(C^*\)-algebra approach \cite{bel17}, employing Gelfand representation. In particular, the work by Heunen et al. \cite{hee0}, \cite{hee} introduces a topos-theoretic framework that redefines the Gelfand spectrum for noncommutative C*-algebras and leading to an explicit Gelfand transform. By constructing a sheaf topos over the poset of commutative subalgebras, their approach introduces an internal Gelfand spectrum that integrates noncommutative algebraic structures with a generalized notion of space. 
This innovative perspective not only provides a bridge between noncommutative and commutative theories but also leads to the formulation of an external spectrum, enabling a geometrical interpretation of noncommutative spaces analogous to classical manifolds.
Selecting an appropriate \(C^*\) algebra \(\mathcal{C}\) for pseudodifferential operators' symbols on manifolds \(\mathcal{M}\), we can choose the spectrum of \(\mathcal{C}\) as the coordinatization method of the manifold, similarly to \cite{bel17}. For this purpose, \cite{hee} not only adapts the concept of Dirac measures to noncommutative settings but also showcases how the algebra's genericity, i.e. its (topologized) spectrum being homeomorphic to the manifold, emerges from the seamless topological integration of \(\mathcal{M}\) with the external spectrum of the algebra.

Hence, in this paper, we prove the following theorems, whose proofs are the topics of Sections \ref{proof1} and \ref{proof2}, respectively:

\begin{thm} \label{ell thm}
The BC method can be used to directly prove uniqueness arguments for inverse elliptic BVPs, when the corresponding edge problem is elliptic in the sense of Schulze’s pseudo-differential operator edge calculus (i.e. both the principal interior symbol and the edge symbol are invertible). 
\end{thm}

\begin{thm} \label{par thm}
The BC method can be used to directly prove uniqueness arguments for inverse Volterra parabolic BVPs, when the corresponding edge problem is elliptic in the sense of Schulze’s pseudo-differential operator edge calculus (i.e. both the principal interior symbol and the edge symbol are invertible).
\end{thm}

Led by the insights provided by the previous two theorems, we propose the following conjecture as a topic worthy of further investigation:

\begin{cl} \label{main thm}
The BC method potentially unifies inverse BVPs of all types, namely hyperbolic, elliptic, and parabolic.    
\end{cl} 

Finally, based on \ref{ell thm}, we demonstrate that the Calderón problem can be studied directly via the Boundary Control method without reducing the elliptic problem to a hyperbolized problem \cite{DosS}.

We now proceed to state the Calderón problem \cite{bal} for a compact smooth Riemannian manifold $\mathcal{M}$ of arbitrary dimension $n \geq 2$.
The direct problem is as follows:
\begin{equation} \label{dir}
\begin{array}{ll}
\nabla \cdot \sigma(x) \nabla u(x) = 0, & x \in \mathcal{M}\\
u(x) = f(x), & x \in \partial \mathcal{M},
\end{array}
\end{equation}
where $\mathcal{M}$ is a compact smooth Riemannian manifold of dimension $n \geq 2$ with smooth boundary $\partial \mathcal{M}$. The function $\sigma(x)$ denotes the smoothly varying conductivity coefficient on $\mathcal{M}$, and $f(x)$ represents the given Dirichlet boundary condition.
For a smooth function $f$, it has a unique smooth classical solution $u = u^{f}(x)$. 

The transformation from Dirichlet to Neumann data, also known as the voltage-to-current mapping, is described by:
\begin{equation} \label{DN}
\left\{ 
\begin{array}{l}
\Lambda: H^{\frac{1}{2}}(\partial X) \to H^{-\frac{1}{2}}(\partial \mathcal{M}), \\
f(x) \mapsto \Lambda_\sigma[f](x) = \sigma(x) \frac{\partial u}{\partial \nu}(x), 
\end{array}
\right.
\end{equation}
where $(\cdot)_{\nu}$ denotes the outward normal derivative.
The Calderón problem involves deducing the conductivity coefficient, $\sigma$, by utilizing the data from the Dirichlet-to-Neumann map.

The uniqueness argument for the inverse problem is defined as follows:
\begin{thm} \label{Calderon}
Let $M_1$ and $M_2$ be two $n$-dimensional Riemannian manifolds with boundary that are locally flat near the boundary.
Assume $\partial M_1 = \partial M_2$ and that they share the same boundary data, i.e. the same Dirichlet-to-Neumann map (\ref{DN}). 
On each manifold consider the edge-degenerate problem in the sense of Schulze corresponding to (\ref{dir}), and suppose $\sigma \neq 0$. Satisfying the ellipticity conditions of Schulze’s PDO edge calculus (meaning the principal interior symbol and principal edge symbol are invertible), 
$M_1$ and $M_2$ are isometrically diffeomorphic for $\gamma\ne \frac{1}{2}, \frac{3}{2}$, where $\gamma$ is the weight factor in the spaces of the edge-degenerate operators.
\end{thm}

\section{Proofs on BC method extensions} \label{proof}

\subsection{Proof on BC for elliptic problems} \label{proof1}
\subsubsection{Facts to be used in the proof}
All definitions/facts included in the following itemized lists are from \cite{SchEl}.\\

\noindent\textit{Graded Algebra and Principal Symbol}

\begin{itemize}
\item[\text{[D0]}]For a manifold \(\mathcal{M}\) with a \(C^\infty\) structure densely defined on a subset \(\mathcal{M}^\circ\), the algebra \(Diff(\mathcal{M}^\circ)\) includes all differential operators with smooth coefficients on \(\mathcal{M}^\circ\). A subalgebra \(\mathcal{D} \subset Diff(\mathcal{M}^\circ)\) is defined by selecting operators based on their behavior near the singularity set \(\mathcal{M} \setminus \mathcal{M}^\circ\), emphasizing the algebra's responsiveness to the manifold's singular and non-singular regions.\label{item:D0}

\item[\text{[D1]}] We define \( \mathcal{D} \) with an order-based filtration: \(\{0\} = \mathcal{D}_{-1} \subset \mathcal{D}_0 \subset \mathcal{D}_1 \subset \cdots\). This leads to the formation of the associated graded algebra \(\mathrm{gr}\ \mathcal{D} = \bigoplus_{j=0}^{\infty} \mathcal{D}_{j}/\mathcal{D}_{j-1}\), where \(\mathrm{gr}_k\mathcal{D}\) isolates symbols of operators precisely of order \(k\), thus segregating different order effects.\label{item:D1}

\item[\text{[D2]}] The principal symbol of an operator \(D\) in \(\mathcal{D}_k\), noted as \(\Sigma_k(D)\), abstracts the operator's highest-order impact by projecting \(D\) to \(\mathrm{gr}_k\mathcal{D}\) through the map \(\Sigma_k: \mathcal{D}_k \rightarrow \mathrm{gr}_k\mathcal{D}\). This operation discards lower-order influences, highlighting the core action of \(D\).\label{item:D2}

\item[\text{[D3]}] Globally, utilizing a partition of unity to "glue" these local correspondences together, while carefully selecting invertible pseudodifferential operators to circumvent topological obstructions, allows for the extension of local algebraic and geometric properties of \( \mathcal{D} \) across the entire manifold \( \mathcal{M} \). \label{item:D3}
\end{itemize}

\vspace{0.7\baselineskip}
\noindent\textit{Localization of Principal Symbol and Monomorphism}

\begin{itemize}
\item[\text{[L0]}] Let the stretched manifold $M$ be a compactification and extension of the singular space $\mathcal{M}$, 
which "stretches" over the singular points of $\mathcal{M}$, effectively smoothing out the singularities. The subspace \(C^{\infty}(\mathcal{M}) \subset C^{\infty}(M)\) acts as a natural subalgebra within \(\mathcal{D}\), with the commutation relation \([C^{\infty}(\mathcal{M}), \mathcal{D}_{k}] \subset \mathcal{D}_{k-1}\) indicating that \(C^{\infty}(\mathcal{M})\) effectively operates on the graded algebra \(\mathrm{gr}\mathcal{D}\), rendering each \(\mathrm{gr}_{k}\mathcal{D}\) a \(C^{\infty}(\mathcal{M})\)-module, where the left and right actions coincide.\label{item:L0}

\item[\text{[L1]}] Localization at a point \(x \in \mathcal{M}\) targets the action of \( \mathcal{D} \)'s operators at \(x\). This process uses the maximal ideal \(I_x \subset C^{\infty}(\mathcal{M})\), comprised of functions zero at \(x\), to adapt the algebra to local considerations.\label{item:L1}

\item[\text{[L2]}] Through localization, each operator \(D\)'s principal symbol, \(\Sigma_k(D)\), is represented at \(x\) by \(\sigma_x(D)\), determined via \(\Sigma_{kx} = \mathrm{gr}_k\mathcal{D}/I_{x}\mathrm{gr}_k\mathcal{D}\). This representation captures the operator's highest-order behavior near \(x\), effectively tying algebraic properties to manifold's topology.\label{item:L2}

\item[\text{[M1]}] The monomorphism induced by mapping \(\pi: \mathrm{gr}_k \mathcal{D} \rightarrow \prod_{x \in \mathcal{M}} \Sigma_{kx}\), with each \(\pi_x: \mathrm{gr}_k \mathcal{D} \rightarrow \Sigma_{kx}\) acting as a natural projection, ensures that the set of local representations \(\{\sigma_x(D)\}\) uniquely identifies the principal symbol of \(D\), thereby establishing a direct link between the algebra's elements and their manifestations across \(\mathcal{M}\).\label{item:M1}
\end{itemize}

\vspace{0.7\baselineskip}
\noindent\textit{Scaling for Localized Behavior}

\begin{itemize}
\item[\text{[S1]}] To discern the influence of \(D\) at \(x\), scaling transformations \(g_{\lambda}\) for \(\lambda \in \mathbb{R}_+\) adjust \(D\) for local behavior analysis. The limit \(D_{x} = \lim_{\lambda \rightarrow \infty} \lambda^{-m}(g_{\lambda}^{*})^{-1}Dg_{\lambda}^{*}u\), for \(u \in C_{0}^{\infty}(K_{x})\) and excluding \(x\) for singular points, showcases \(D\)'s localized impact, correlating with \(\sigma_x(D)\).\label{item:S1}
\end{itemize}

\vspace{0.7\baselineskip}
\noindent\textit{Fundamental Relationship}

\begin{itemize}
\item[\text{[T1]}] The relationship \(\sigma_x(D) \mapsto D_{x}\) is both 1-1 and multiplicative, highlighting the principal symbol's role in unveiling \(D\)'s essence in a localized context. This interplay, facilitated by the monomorphism in [M1] and the scaling process in [S1], deepens our comprehension of differential operators' localized effects, affirming the algebra \(\mathcal{D}\)'s adaptability and the continuity of its application across \(\mathcal{M}\).\label{item:T1}
\end{itemize}

\vspace{0.7\baselineskip}
\noindent\textit{Edge Symbols}

\begin{itemize}
\item[\text{[E1]}] Edge symbols, represented as \( D(x; \xi) \), act on the specific weighted Sobolev spaces \( K^s_{\gamma}(K_{\Omega}) \), which are constructed by gluing weighted and standard Sobolev spaces across different regions of the cone \( K_{\Omega} \).
\( K^s_{\gamma}(K_{\Omega}) \) is an adaptation for cone structures of the abstract wedge spaces \(\mathcal{W}^s(\mathcal{M}; K_{\Omega})\).
The connection between the algebraic structure of edge symbols and  \(\mathcal{W}^s(\mathcal{M}; K_{\Omega})\)
underlies the continuity and scaling properties essential for the functionality of edge-degenerate pseudodifferential operators within the \( \mathcal{W}^s(\mathcal{M}; K_{\Omega}) \) spaces.
\label{item:E1}
\end{itemize}

\vspace{0.7\baselineskip}
\noindent\textit{Quantization}

\begin{itemize}
    \item[\text{[Q1]}] The exact sequence 
    \[0 \rightarrow J_{\gamma}^{m} \rightarrow \mathcal{A}_{\gamma}^{m} \rightarrow \mathcal{O}^{m} \rightarrow 0\]
    delineates the structural relationships in the algebra of quantization for pseudodifferential operators. Here, \(J_{\gamma}^{m}\) is an ideal comprising edge symbols, \(\mathcal{A}_{\gamma}^{m}\) is an algebra of both interior and edge symbols meeting specific compatibility conditions, and \(\mathcal{O}^{m}\) consists of interior symbols, enabling a comprehensive formulation of quantization processes in the context of manifold \(\mathcal{M}\).\label{item:Q1}
\end{itemize}

\vspace{0.7\baselineskip}
\noindent\textit{Calkin Algebra}

\begin{itemize}
\item[\text{[A1]}]  Let \(\mathcal{H}\) denote a Hilbert space, and \(\mathcal{B}(\mathcal{H})\) the algebra of bounded operators on \(\mathcal{H}\). The set \(\mathcal{\tilde{A}} \subset \mathcal{B}(\mathcal{H})\), representing zero-order pseudodifferential operators (\(\psi\mathrm{D}\mathrm{O}\)), forms a subalgebra known as the general algebra of zero-order \(\psi\mathrm{D}\mathrm{O}\)s (the general case follows
by order reduction).\label{item:A1}

\item[\text{[A2]}] A homomorphism \(\sigma: \mathcal{\tilde{A}} \rightarrow S\) is defined, mapping \(\mathcal{\tilde{A}}\) into a unital topological algebra \(S\). This mapping is characterized by \(\sigma(A) = 0\) if and only if \(A \in \mathcal{K}\), where \(\mathcal{K}\) denotes the ideal of compact operators within \(\mathcal{B}(\mathcal{H})\).\label{item:A2}

\item[\text{[A3]}] The symbol mapping \(\sigma\) induces a well-defined monomorphism \(\tilde{\sigma}: \mathcal{\tilde{A}}/(\mathcal{K} \cap \mathcal{\tilde{A}}) \rightarrow S\), effectively establishing a one-to-one correspondence between the elements of the Calkin algebra \(\mathcal{\tilde{A}}/(\mathcal{K} \cap \mathcal{\tilde{A}})\) and the elements within \(S\).\label{item:A3}

\item[\text{[A4]}] There exists a continuous linear mapping \(Q: S \rightarrow \mathcal{\tilde{A}}\), for which \(\sigma(Q(s)) = s\) for all \(s \in S\), making \(Q\) the right inverse of \(\sigma\). This mapping is instrumental in ensuring that the structure and relational properties between \(\mathcal{\tilde{A}}\) and \(S\) are preserved.\label{item:A4}
\end{itemize}

\vspace{0.7\baselineskip}
\noindent\textit{Edge Morphisms}

\begin{itemize}
    \item[\text{[EM1]}] A morphism \(\mathrm{A}\) of order \(m\) and weight \(\gamma\), represented as an operator matrix \(\mathrm{A} = \left( \begin{smallmatrix} A & C \\ B & D \end{smallmatrix} \right)\), includes \(A\) as a pseudodifferential operator on \(\mathcal{M}\), with \(B\) and \(C\) as edge boundary and coboundary operators, and \(D\) as a pseudodifferential operator on the edge \(X\), all of order \(m\) and weight \(\gamma\). The set of such morphisms is denoted by \(\mathrm{Mor}_\gamma^m(\mathcal{M})\). The interior symbol of morphism \(\mathrm{A}\), \(\sigma(\mathrm{A})\), is defined as \(\sigma(A)\), the interior symbol of the operator \(A\).\label{item:EM1}
\end{itemize}

\subsubsection{Lemmas}
In the framework of a noncommutative Calkin algebra on a manifold, we deploy a topos-theoretic approach in accordance to \cite{hee}.
The internal spectrum models a C*-algebra's structure within a sheaf topos, and the external spectrum effectively translates these algebraic points into a topological structure that can be visualized as points in a space.
It is inferred from \cite{hee} that norm convergence in an algebra $A$ implies topological convergence in the external spectrum \(\Sigma_{\text{ext}} A\), as the latter carries the weak* topology, which respects norm limits. In particular, the topology of the internal spectrum, which is aligned with weak* features suitable for C*-algebra structures, is reflected in the external spectrum to maintain mathematical consistency and meaningful interpretation.
Based on this inference the map \(\epsilon: \mathcal{M} \to \Sigma_{\text{ext}}A\), where \(\epsilon_x=\delta_x\) represents the Dirac measure is meaningful.

The Calkin algebra \(\mathcal{C}\) on \(\mathcal{M}\), defined in \hyperref[item:A3]{[A3]}, achieves a notion of 'genericity' through the homeomorphic relation established between \(\mathcal{M}\) and \(\Sigma_{ext}(\mathcal{C})\) via the embedding \(\epsilon\), enriched by 
a topos-theoretic generalization of spatial isomorphism; the Gelfand transform \(\hat{\mathcal{C}} \cong \mathcal{C}\) in this context reflects a deep algebraic-geometric correspondence, emphasizing the continuity and topological equivalence between \(\mathcal{M}\) and \(\Sigma_{ext}(\mathcal{C})\).

Additionally in \cite{hee} the broadened perspective on spectrum and space through sheaf toposes and frames offers the flexibility to rigorously define and use the external spectrum for non-C*-algebraic structures, if it can be described within a suitable framework (like a topos or through frame theory).
Hence, the algebra \(S\) (\hyperref[item:A2]{[A2]}, \hyperref[item:A3]{[A3]}) is eligible for external spectrum. In particular, defining \(\epsilon: M \rightarrow \Sigma_{ext}(\mathcal{C})\), where \(\epsilon_x=\delta_x\) is the Dirac measure, the continuity of \(\epsilon\) and its inverse suggest that \(S\) respects the spatial topology required for an external spectrum; this continuity ensures that \(S\) can be effectively treated as a fibered space over the manifold, akin to the external spectrum's behavior over the base space of subalgebras.
Also the algebraic operations within \(S\) need to be compatible with the fibered structure; the monomorphism \hyperref[item:A3]{[A3]} and symbol mapping \hyperref[item:A2]{[A2]} ensure that \(S\) retains algebraic properties consistent with the Calkin algebra, crucial for its interpretation as an external spectrum.
Of course, an internal spectrum cannot be constructed for \(S\) due to its potentially non-C*-algebra properties, hence S has no explicit Gelfand transform.

Let us note that S in \hyperref[item:A2]{[A2]} can refer to either interior or edge symbol algebras of the edge-degenerate PDOs. However, in what follows, for ease of writing, by abuse of terminology, we will refer to $S$ as the algebra of principal symbols (standard ones; not the ones for edge-degenerate operators), instead of the algebra of interior symbols, as the two essentially coincide. 

\begin{lem}  \label{11S}
There is a one-to-one correspondence between points on the manifold \(\mathcal{M}\) and the external spectrum \(\Sigma_{ext}(S)\) of principal symbols \(S\), which is related to differential operators on \(\mathcal{M}\).
\end{lem}
\noindent \textit{Proof}.

\textit{Local Association Through Localization.} (\hyperref[item:L1]{[L1]}, \hyperref[item:L2]{[L2]}). 
For any point \(x\) on \(\mathcal{M}\), the process of localization maps the behavior of differential operators at \(x\). This involves identifying the principal symbols \(\sigma_x(D)\) for each differential operator \(D\) at \(x\), using maximal ideals \(I_x\) comprised of functions that are zero at \(x\). 

\textit{Global Extension and Correspondence.}
The monomorphism \hyperref[item:M1]{[M1]} expands these local associations to a global structure within \(S\). It guarantees that the connections between points \(x\) on \(\mathcal{M}\) and their localized algebraic behaviors via \(\sigma_x(D)\) correspond to unique maximal ideals in \(S\) without overlap. This ensures a one-to-one correspondence between points on \(\mathcal{M}\) and maximal ideals in \(S\), summarizing the differential actions at each point \(x\). 
\(\hfill \square \)

\begin{lem} \label{11Calkin}
There is one-to-one correspondence between the manifold \(\mathcal{M}\) and 
the external spectrum \(\Sigma_{ext}(\mathcal{C})\) of the Calkin algebra, \(\mathcal{\tilde{A}}/(\mathcal{K} \cap \mathcal{\tilde{A}})\).
\end{lem}
\noindent \textit{Proof}. By \hyperref[item:A3]{[A3]}, \(\tilde{\sigma}\) creates a one-to-one correspondence between maximal ideals in the Calkin algebra and \(S\), preserving their structural and relational properties. By Lemma \ref{11S} 
the external spectrum \(\Sigma_{ext}(S)\) of \(S\) and the manifold are in 1-1 correspondence.
Therefore, there is one-to-one correspondence between the manifold and the external spectrum of the Calkin algebra. $\hfill \square$

\begin{lem} \label{conteS}
The map \(\epsilon: \mathcal{M} \to \Sigma_{ext} S\), where \(\epsilon_x=\delta_x\) is the Dirac measure, associating each point \(x\) in the manifold \(\mathcal{M}\) with a maximal ideal of the external spectrum \(\Sigma_{ext} S\) of the algebra \( S \), is continuous.
\end{lem}
\noindent \textit{Proof}. In accordance to the constructions in \cite{hee}, norm convergence in $S$ implies topological convergence in the external spectrum \(\Sigma_{\text{ext}} S\); as was explained in the beginning of this section.

\textit{Evaluation Map Definition for \( S \)}.  
Define \( ev_x(\sigma(D)) \) for any \( x \in \mathcal{M} \) and a symbol \(\sigma(D) \in S\) to denote the action of the differential operator \( D \) at \( x \), modulo lower order terms, reflecting the graded nature of \( S \).

\textit{Density and Approximation Argument}.  
Given \(\sigma(D) \in S\) and any \(\epsilon > 0\), we utilize the density of \(\mathcal{D}\) in \(\tilde{\mathcal{A}}\) \hyperref[item:D0]{[D0]} to select \(D' \in \mathcal{D}\) such that \(\sigma(D')\) approximates \(\sigma(D)\) in the norm of \( S \), i.e., \(\|\sigma(D) - \sigma(D')\|_S < \epsilon/3\).

\textit{Formal Convergence Argument}.  
For a converging sequence \(x_n \rightarrow x\) in \(\mathcal{M}\), and considering \(\sigma(D') \in S\) that closely approximates \(\sigma(D)\) in \( S \):
\begin{eqnarray*}
\left\| \sigma(D)(x_n) - \sigma(D)(x) \right\|_S &\leq & \left\| \sigma(D)(x_n) - \sigma(D')(x_n) \right\|_S + \left\| \sigma(D')(x_n) - \sigma(D')(x) \right\|_S \\
&& + \left\| \sigma(D')(x) - \sigma(D)(x) \right\|_S,
\end{eqnarray*}
where the continuity of \(D'\) ensures \(\left\| \sigma(D')(x_n) - \sigma(D')(x) \right\|_S \rightarrow 0\) as \(n \rightarrow \infty\).

\textit{Ensuring Continuity}.  
This formulation guarantees that \(\left\| \sigma(D)(x_n) - \sigma(D)(x) \right\|_S < \epsilon\), affirming the continuity of \(\epsilon\). 
$\hfill \square$

\begin{lem} \label{conte}
The map \(\epsilon: \mathcal{M} \to \Sigma_{ext}\mathcal{C}\), where \(\epsilon_x=\delta_x\) is the Dirac measure, associating each point \(x\) in the manifold \(\mathcal{M}\) 
with a maximal ideal of the external spectrum \(\Sigma_{ext}\mathcal{C}\) of the Calkin algebra \(\mathcal{C}\) is continuous.
\end{lem}
\noindent \textit{Proof}.
The argument is the same with the one in Lemma \ref{conteS} applied to \([A] \in \mathcal{C}\), denoting the action of \(A\) modulo compact operators at \(x\) in accordance to the algebra’s quotient structure, instead of \(\sigma(D) \in S\), and with the density argument:

Given \([A] \in \mathcal{C}\) and any \(\epsilon > 0\), we utilize the density of \(\mathcal{D}\) in \(\tilde{\mathcal{A}}\) to select \(D \in \mathcal{D}\) such that \([D]\) approximates \([A]\) in the quotient norm of \(\mathcal{C}\), i.e., \(\|[A] - [D]\|_{\mathcal{C}} < \epsilon/3\). $\hfill \square$

\begin{lem} \label{conte1S}
The map \(\epsilon^{-1}: \Sigma_{ext} S \rightarrow \mathcal{M}\), where \(\epsilon_x=\delta_x\) is the Dirac measure, from a maximal ideal of the external spectrum \(\Sigma_{ext} S\) of the algebra \( S \), to a point in the manifold \(\mathcal{M}\), is continuous.
\end{lem}
\noindent \textit{Proof}. In accordance to the constructions in \cite{hee}, norm convergence in $S$ implies topological convergence in the external spectrum \(\Sigma_{\text{ext}} S\); as was explained in the beginning of this section.

\textit{Sequence of Ideals and Convergence}.  
Assume a sequence of maximal ideals \(M_n\) converges to \(M\) in \(\Sigma_{ext} S\). This implies that for every principal symbol \(\sigma(D) \in S\), the evaluation at \(M_n\) (which can be thought of as evaluating \(\sigma(D)\) at the points \(x_n\) corresponding to \(M_n\)) converges to the evaluation at \(M\) (evaluation at the point \(x\) corresponding to \(M\)) as \(n \rightarrow \infty\), i.e.
\[
\sigma(D)|_{M_n} \rightarrow \sigma(D)|_M \text{ as } n \rightarrow \infty.
\]

\textit{Approximation and Density in \(\mathcal{D}\)}.  
For any \(\sigma(D) \in S\) and given \(\epsilon > 0\), there exists an operator \(D' \in \mathcal{D}\) approximating \(\sigma(D)\) within \(S\), such that \(\|\sigma(D) - \sigma(D')\|_S < \epsilon/3\). This uses the density of \(\mathcal{D}\) in \(\tilde{\mathcal{A}}\)
\hyperref[item:D0]{[D0]}. 

\textit{Continuity of \(\sigma(D')\) and Convergence of Points}.  
Consider \(x_n, x \in \mathcal{M}\) corresponding to \(M_n, M\), respectively. The continuity of \(D'\) implies \(D'(x_n) \rightarrow D'(x)\). Since \(\sigma(D')\) reflects \(D'\)'s highest-order action, this continuity extends to \(\sigma(D')\), ensuring that the evaluations \(\sigma(D')(x_n)\) converge to \(\sigma(D')(x)\) within the algebraic framework of \(S\).

\textit{Ensuring Continuity of \(\epsilon^{-1}\)}.  
The convergence \(\sigma(D')(x_n) \rightarrow \sigma(D')(x)\), due to the operational continuity of \(D'\) and its approximation of \(\sigma(D)\), ensures that the spectral convergence \(M_n \rightarrow M\) mirrors as the topological convergence \(x_n \rightarrow x\) on \(\mathcal{M}\), thereby demonstrating the continuity of \(\epsilon^{-1}\). $\hfill \square$

\begin{lem} \label{conte1}
The map \(\epsilon^{-1}: \Sigma_{ext}\mathcal{C} \rightarrow \mathcal{M}\), where \(\epsilon_x=\delta_x\) is the Dirac measure, from a maximal ideal of the external spectrum \(\Sigma_{ext}\mathcal{C}\) of the Calkin algebra \(\mathcal{C}\), to a point in the manifold \(\mathcal{M}\) is continuous.
\end{lem}
\noindent \textit{Proof}.
The argument is the same with the one in Lemma \ref{conte1S} applied to \([f] \in \mathcal{C}\) instead of \(\sigma(D) \in S\), and with the density argument:

For any \([f] \in \mathcal{C}\) and given \(\epsilon > 0\), there exists an operator \(D \in \mathcal{D}\) approximating \([f]\) within \(\mathcal{C}\), such that \(\|[f] - [D]\|_{\mathcal{C}} < \epsilon/3\). This utilizes the density of \(\mathcal{D}\) in \(\tilde{\mathcal{A}}\) \hyperref[item:D0]{[D0]}. 
$\hfill \square$

\begin{lem} \label{comS}
The commutative structure of the differential operators algebra \( \mathcal{D} \) is preserved in the algebra of principal symbols \( S \).
\end{lem}
\noindent \textit{Proof}. By contradiction, assume the commutative structure of \( \mathcal{D} \) is not preserved in \( S \). This would imply a failure in the commutative properties' transition from \( \mathcal{D} \) to \( S \).

\textit{Commutative Structure of \(\mathcal{D}\)}.
The commutative properties of \(\mathcal{D}\) with respect to operations involving \(C^{\infty}(\mathcal{M})\) are detailed in \hyperref[item:L0]{[L0]}. The overall architecture of \(\mathcal{D}\) mirrors these local commutative behaviors, supported by an order-based filtration (\hyperref[item:D1]{[D1]}) and reinforced by employing a partition of unity along with the choice of invertible pseudodifferential operators to avoid topological obstructions (\hyperref[item:D3]{[D3]}).

\textit{Preservation of Commutativity in \(S\)}.
Through the monomorphism \hyperref[item:M1]{[M1]}, local principal symbols \(\sigma_x(D)\) are systematically connected to \(S\), preserving local commutative structures globally. The relationship between \(\sigma_x(D)\) and \(D_x\) (\hyperref[item:T1]{[T1]}) extends commutative properties from local to global contexts within \( S \).

\textit{Contradiction and Conclusion}.  
The mechanisms ensuring the preservation of commutativity from \( \mathcal{D} \) through \( S \), as outlined by the structured mapping [M1], contradict the assumption of disrupted commutative structures. Therefore, the commutative structure within \( \mathcal{D} \) is indeed preserved in the algebra of principal symbols \( S \). $\hfill \square$

\begin{lem} \label{comCalkin}
The commutative structure of the differential operators algebra \(\mathcal{D}\) is preserved in the Calkin algebra \(\mathcal{\tilde{A}}/(\mathcal{K} \cap \mathcal{\tilde{A}})\).
\end{lem}
\noindent \textit{Proof}. By contradiction,
assume the commutative structure of \(\mathcal{D}\) is not preserved in the Calkin algebra. 
However, by Lemma \ref{comS} the commutative structure within \( \mathcal{D} \) is preserved in the algebra of principal symbols \(S\).
Hence, the mechanism ensuring commutativity's preservation from \(\mathcal{D}\) through \(S\), and by \hyperref[item:A3]{[A3]}, to the Calkin algebra, contradict the assumption of disrupted commutative structures. Therefore, the commutative structure within \(\mathcal{D}\) is indeed preserved in the Calkin algebra, \(\mathcal{\tilde{A}}/(\mathcal{K} \cap \mathcal{\tilde{A}})\). $\hfill \square$

\begin{lem} \label{ACCS}
The algebra of principal symbols \(S\), associated with the differential operators on the manifold \( \mathcal{M} \), satisfies the Ascending Chain Condition (ACC).
\end{lem}
\noindent \textit{Proof}. By contradiction,
assume that \(S\) does not satisfy the ACC. This implies the existence of an infinite ascending sequence \(\{I_n\}\) of ideals in \(S\) that does not stabilize, i.e. \(I_n \subset I_{n+1}\) for all \(n\) and \(I_n \neq I_{n+1}\).

\textit{Local Boundedness of \(D_x\)}. 
According to \hyperref[item:S1]{[S1]}, for each \(x\) in \(\mathcal{M}\), the operational impact of \(D_x\) is bounded due to scaling transformations, indicating that differential operators exhibit locally bounded complexity.

\textit{Transition of Boundedness from \(D_x\) to \(\sigma_x(D)\)}. 
The monomorphism \hyperref[item:T1]{[T1]} illustrates that each local principal symbol \(\sigma_x(D)\) directly correlates to the bounded limit \(D_x\), implying that the complexity of \(\sigma_x(D)\) is similarly bounded at each point on \(\mathcal{M}\).

\textit{Connection of Local Symbols \(\sigma_x(D)\) to Global \(S\)}. Through the monomorphism \hyperref[item:M1]{[M1]}, local principal symbols \(\sigma_x(D)\) are systematically linked to their global counterparts within \(S\), ensuring the transfer of local bounded complexity to the algebra of principal symbols \(S\).

\textit{Implication of Filtration on Global \(D\) and \(S\)}. 
The filtration of \(\mathcal{D}\) as detailed in \hyperref[item:D1]{[D1]} not only organizes differential operators by their order but also enables the translation of local boundedness observed through \(D_x\) and \(\sigma_x(D)\) to a structured complexity within the global framework of \(\mathcal{D}\) and subsequently to \(S\). 
This global structure of \(\mathcal{D}\)'s structure is further ensured through the utilization of a partition of unity and circumventing topological obstruction by selection of invertible pseudodifferential operators (\hyperref[item:D3]{[D3]}).

\textit{Contradiction and Conclusion}.
The inherent bounded complexity of symbols in $S$ implies that any ascending chain of ideals in \(S\) must stabilize, contradicting the initial assumption. Therefore, \(S\) must satisfy ACC. $\hfill \square$

\begin{lem} \label{ACCCalkin}
The Calkin algebra \( \mathcal{C}=\mathcal{\tilde{A}}/(\mathcal{K} \cap \mathcal{\tilde{A}})\) satisfies the Ascending Chain Condition (ACC).
\end{lem}
\noindent \textit{Proof}. By contradiction, assume that the Calkin algebra $\mathcal{C}$ does not satisfy the ACC, leading to an infinite ascending chain of ideals \(\{I_n\}\).

\textit{Monomorphism and ACC in \(S\)}.
By \hyperref[item:A3]{[A3]}, \(\tilde{\sigma}\) creates a one-to-one correspondence between ideals in the Calkin algebra and \(S\), preserving their structural and relational properties. Since \(S\) satisfies the ACC by Lemma \ref{ACCS}, every ascending chain of ideals stabilizes.

\textit{Contradiction and Conclusion}. 
The infinite ascending chain \(\{I_n\}\) in \(\mathcal{\tilde{A}}/(\mathcal{K} \cap \mathcal{\tilde{A}})\) would correspond to a similar chain in \(S\) through \(\tilde{\sigma}\), contradicting \(S\)'s ACC property. The contradiction confirms our initial assumption is false; therefore, the Calkin algebra must satisfy the ACC. $\hfill \square$

\begin{lem} \label{genericC}
The Calkin algebra $\mathcal{C}$ is a generic algebra, i.e. the embedding \(\epsilon: \mathcal{M} \to \Sigma_{ext}\mathcal{C}\), where \(\epsilon_x=\delta_x\) is the Dirac measure, exists and is a homeomorphism, and the Gelfand transform of $\mathcal{C}$ can be given explicitly.
\end{lem}
\noindent \textit{Proof}. 
The Calkin algebra is eligible for application of the paper \cite{hee} on construction of Gelfand transform for noncommutative $C^*$ algebras, as it contains commutative subalgebras by Lemma \ref{comCalkin} and it satisfies the Ascending Chain Condition by Lemma \ref{ACCCalkin}.

The one-to-one correspondence between points on \( \mathcal{M} \) and maximal ideals in the external spectrum of the Calkin algebra \(\mathcal{C}\) is given by Lemma \ref{11Calkin}.
The continuity of the embedding \( \epsilon \) and its inverse \( \epsilon^{-1} \) is given by Lemmas \ref{conte} and \ref{conte1}, respectively. 

The link between external and internal spectrum as given in \cite{hee} allows for the explicit construction of Gelfand transform. $\hfill \square$

\begin{lem} \label{genericS}
The algebra $S$ of principal symbols is a generic algebra, i.e. the embedding \(\epsilon: \mathcal{M} \to \Sigma_{ext}S\), where \(\epsilon_x=\delta_x\) is the Dirac measure, exists and is a homeomorphism.
\end{lem}
\noindent \textit{Proof}. 
The one-to-one correspondence between points on \( \mathcal{M} \) and maximal ideals in the external spectrum of algebra \(S\) are given by Lemma \ref{11S}.

The continuity of the embedding \( \epsilon \) and its inverse \( \epsilon^{-1} \) are given by Lemmas \ref{conteS} and \ref{conte1S}, respectively. $\hfill \square$

\begin{lem} \label{coord}
The coordinatization of the manifold can be given via the Gelfand transform of the Calkin algebra.
\end{lem}
\noindent \textit{Proof}. 
By Lemmata \ref{genericC} and \ref{genericS} the manifold is homeomorphic to the external spectrum of the Calkin algebra and the external spectrum of the algebra of principal symbols $S$, respectively. Hence, the external spectrum of the Calkin algebra and the external spectrum of $S$ are also homeomorphic. Given that, as $S$ is not eligible for an internal spectrum, whereas the Calkin algebra can support both internal and external spectral analysis, we will use the Calkin algebra for coordinatizing the manifold by its Gelfand transform. $\hfill \square$

\begin{lem} \label{isis}
Given two exact and splitting sequences
\begin{eqnarray}
0 \rightarrow J1_\gamma^m \rightarrow A1_\gamma^m \rightarrow O1^m \rightarrow 0,\\
0 \rightarrow J2_\gamma^m \rightarrow A2_\gamma^m \rightarrow O2^m \rightarrow 0,
\end{eqnarray}
and the isometric isomorphism $\phi: J1_\gamma^m \rightarrow J2_\gamma^m$, then $O1^m$ is isometrically isomorphic to $O2^m$. 
\end{lem}
\noindent \textit{Proof}:
We can express $A1_\gamma^m$ and $A2_\gamma^m$ as direct sums due to the exactness and splitting of the sequences:
\begin{eqnarray}
A1_\gamma^m = J1_\gamma^m \oplus O1^m, \\
A2_\gamma^m = J2_\gamma^m \oplus O2^m. 
\end{eqnarray}

Define a mapping $\psi: A1_\gamma^m \rightarrow A2_\gamma^m$ where $\psi(j, o) = (\phi(j), o)$ for any $(j, o) \in A1_\gamma^m$. This map $\psi$ is constructed to maintain the structural relationship between $A1_\gamma^m$ and $A2_\gamma^m$.

To verify $\psi$ as an isometric isomorphism, we check for linearity, injectivity, surjectivity, and isometry:

1. Linearity: $\psi$ respects vector addition and scalar multiplication, given its operation on each component of the direct sum.

2. Injectivity: Assume $\psi(j_1, o_1) = \psi(j_2, o_2)$. This leads to $(\phi(j_1), o_1) = (\phi(j_2), o_2)$. Given $\phi$ is an isomorphism, we deduce $j_1 = j_2$ and thus $(j_1, o_1) = (j_2, o_2)$, proving injectivity.

3. Surjectivity: For any $(j_2, o_2) \in A2_\gamma^m$, there exists $(j_1, o_2) \in A1_\gamma^m$ such that $\psi(j_1, o_2) = (\phi (j1), o2) = (j_2, o_2)$, indicating surjectivity.

4. Isometry: Since \(\phi\) is an isometry between \(J1_\gamma^m\) and \(J2_\gamma^m\), preserving inner product structures, the mapping \(\psi\), defined as \(\psi(j, o) = (\phi(j), o)\), naturally extends this isometry to \(A1_\gamma^m\) and \(A2_\gamma^m\). By applying \(\phi\) to the \(J\) component and directly transferring the \(O\) component, \(\psi\) maintains the inner product relations intact across the entire structure.

Since $\psi$ and $\phi$ are isometric isomorphisms between $Ai_\gamma^m$'s and $Ji_\gamma^m$'s, $i=1,2$, then   $Oi^m=Ai_\gamma^m/Ji_\gamma^m,\ i=1,2$  are isometrically isomorphic. $\hfill \square$ 

\subsubsection{Proof of Theorem \ref{ell thm}} \label{proofell}

We consider two manifolds $\mathcal{M}_i,\ i=1,2$ with the same boundary and same boundary data, i.e. having the same Dirichlet-to-Neumann map. 

\textbf{1.} On the two manifolds $\mathcal{M}_i,\ i=1,2$, let us consider the following direct BVPs $P_i=\left(\begin{array}{c}
A\\
B_i
\end{array}\right),\ \text{for}\ i=1,2$,
where $A$ is our elliptic operator, $B_1$ is the boundary condition with the input data of the inverse problem, i.e the Dirichlet data, and $B_2$ corresponds to the output data, i.e. the Neumann data, connected by the given Dirichlet-to-Neumann map $\Lambda$ which the two manifolds share.

This is meaningful, as in the static context of elliptic inverse BVPs, both the input (Dirichlet) and output (Neumann) boundary data are snapshots of the steady
state, not states in a temporal evolution.
Hence, using the output data as boundary conditions for a conceptual direct problem does not inherently involve reverse engineering in a temporal sense, because the data does not inherently possess a "direction" of propagation; unlike the dynamic, time-evolving problem.  

\textbf{2.} We will treat the elliptic BVPs $P_i,\ i=1,2$, as edge problems, where the edge is the boundary and the cone degenerates
to $\mathbb{R}_+$, interpreted as the inner normal. 
 By \hyperref[item:EM1]{[EM1]} an operator $A$ and its edge morphism have the same interior symbol, i.e. $\sigma(\mathrm{Ai})=\mathcal{O}_i,\ i=1,2\), where $\mathrm{Ai}\in\mathrm{Mor}^0_{\gamma}$, and \(\mathcal{O}_i,\ i=1,2\) denote the interior symbols of the operator $A$ for the corresponding manifold (they are determined by the manifold). 

\textbf{3.} From \cite{PsMSi}, Section 3.1.2, Proposition 3 and the comments following it, we get 

\begin{prop} \label{PsMSi_Proposition 3}
(\cite{PsMSi}) {\it The operator} $Z=F_{\eta\rightarrow y}^{-1}\kappa^{-1}(\eta)F_{y\rightarrow\eta}$, {\it where} $\kappa$ {\it is a continuous groups with the strong operator topology}, {\it defines an isometric isomorphism}
\begin{equation} \label{PsMSi_(7)}
Z:\mathcal{W}^{s}(\mathbb{R}^{\mathrm{q}},\ E)\rightarrow H^{s}(\mathbb{R}^{\mathrm{q}},\ E).   
\end{equation}
\end{prop}
$\mathcal{W}^{s}$ is the abstract wedge space defined in \hyperref[item:E1]{[E1]}, which is equivalent to the boundary in our case. 

Also wedge spaces modelled by $\{E,\ \kappa_{\lambda}\}, \{E,\tilde{\kappa}_{\lambda}\}$ are isomorphic for any choice of continuous groups with the strong operator topology. If
$$
 \tilde{Z}:\mathcal{W}^{s}(\mathbb{R}^{\mathrm{q}},\ E)\rightarrow H^{s}(\mathbb{R}^{\mathrm{q}},\ E)
$$
is the isomorphism of Proposition \ref{PsMSi_Proposition 3} related to $\tilde{\kappa}_{\lambda}$, then we get the isometric isomorphism
$$
 \tilde{Z}^{-1}Z:\mathcal{W}^{s}(\mathbb{R}^{\mathrm{q}},\ E)\rightarrow \tilde{\mathcal{W}}^{s}(\mathbb{R}^{\mathrm{q}},\ E).
$$

Given the isometric isomorphism between the abstract wedge spaces \( W^s(R^q, E) \) and \( \widetilde{W}^s(R^q, E) \) established through this specific isomorphism \( \widetilde{Z}^{-1}Z \), 
the edge symbols \( D(x; \xi) \) on \( K^s_{\gamma}(K_{\Omega}) \) and \( \widetilde{D}(x; \xi) \) on \( \widetilde{K}^s_{\gamma}(K_{\Omega}) \) are isometrically isomorphic. This is because as stated in \hyperref[item:E1]{[E1]} 
\( K^s_{\gamma}(K_{\Omega}) \) is an adaptation of \( W^s(R^q, E) \) with $E=K_{\Omega}$ to specific geometric configurations, i.e. cones, and both sets of spaces serve analogous roles in the context of edge-degenerate operators and their symbol calculus (in particular, these spaces define the domains and codomains for the operators, and their transformations are dictated by similar structural and functional principles).

\textbf{4.} Applying Lemma \ref{isis} for sequences described in \hyperref[item:Q1]{[Q1]}, and using 3. we conclude that the interior symbols of the manifolds are isometrically isomorphic, i.e. there exists an isometric isomorphism \(\Psi: \mathcal{O}_1 \to \mathcal{O}_2\).

\textbf{5.} By Lemma \ref{coord} the coordinatization of the manifold can be given via the Gelfand transform of the Calkin algebra.
In particular, by Lemma \ref{coord}, the identification between \(M_1\) and \(M_2\) via the algebraic isomorphism of their principal symbol algebras induces a homeomorphism \(\Phi: M_1 \to M_2\). Since the isometric isomorphism 
$\Psi$ of Step 4
preserves the metric structure at each point, the pullback of the metric under \(\Phi\) satisfies \(\Phi^* g_2 = g_1\), where \(g_i,\ i=1,2 \) is the Riemannian metric on \( M_i,\ i=1,2 \),
 ensuring that \(\Phi\) is a distance-preserving bijection. The Myers–Steenrod theorem then guarantees that any bijective isometry between smooth Riemannian manifolds is necessarily a smooth diffeomorphism. Consequently, we conclude that \(\Phi\) is an isometric diffeomorphism, i.e. that \(M_1\) and \(M_2\) are isometrically diffeomorphic.
$\hfill \square$

\subsection{Proof on BC for parabolic problems} \label{proof2}

\subsubsection{Facts to be used in the proof}
All definitions/facts included in the following itemized list are from \cite{SchEl}.

\vspace{0.7\baselineskip}
\noindent\textit{Parabolic Symbols and Volterra Operators}

\begin{itemize} 
    \item[\text{[P1]}] A symbol \(a \in S_{V}^{\mu;\ell}(\mathbb{R}^{n}\times \mathbb{R}^{n}\times \mathbb{H};E,\tilde{E})\) is parabolic if it demonstrates parameter-dependent ellipticity within \(S^{\mu;\ell}(\mathbb{R}^{n}\times \mathbb{R}^{n}\times \mathbb{H};E,\tilde{E})\), ensuring crucial analytical properties over \(\mathbb{H}\) necessary for temporal dynamics in parabolic differential equations. Here, \(\mu\) denotes the order, \(\ell\) the degree of homogeneity, and \(E,\ \tilde{E}\) represent Hilbert spaces. \label{item:P1}

    \item[\text{[P2]}] A Volterra operator \(A(\lambda) \in L_{V\ (cl)}^{\mu;\ell}(X;\mathbb{H};E,\ F)\) qualifies as parabolic when it satisfies parameter-dependent ellipticity within \(L_{(cl)}^{\mu;\ell}(X;\mathbb{H};E,\ F)\). This class of operators, emerging from symbols exhibiting the Volterra property, is specialized for equations typified by temporal parameters, aligning with parabolic equation characteristics.\label{item:P2}

    \item[\text{[P3]}] The sequence of principal symbols for Volterra operators is exact and splits as follows:
    \[
    0 \rightarrow L_{V\ cl}^{\mu-1;\ell}(X;\mathbb{H};E,\ F)\stackrel{\iota}{\rightarrow} L_{V\ cl}^{\mu;\ell}(X;\mathbb{H};E,\ F)\stackrel{\sigma_{\psi}^{\mu;\ell}}{\rightarrow} S_{V}^{(\mu;\ell)}((T^{*}X\times \mathbb{H})\backslash 0,\ \mathrm{Hom}(\pi^{*}E,\ \pi^{*}F)) \rightarrow 0
    \]
    Here, \(\sigma_{\psi}^{\mu;\ell}\) denotes the principal symbol, \(\pi^{*}\) refers to the pull-back associated with the projection \(\pi\), and \(S_{V}^{(\mu;\ell)}\) represents the space of anisotropic homogeneous functions of degree \(\mu\) that are analytic within \(\mathbb{H}\). \label{item:P3}
\end{itemize}

\subsubsection{Proof on Theorem \ref{par thm}}
Through \hyperref[item:P1]{[P1]} and \hyperref[item:P2]{[P2]} Volterra operators' parabolicity is linked with parameter-dependent ellipticity, hence we resort to the solution for inverse elliptic BVPs in Section \ref{proofell}. 
Theorem \hyperref[item:P3]{[P3]} is required at step 4 to ensure that symbols are isomorphic to the Calkin algebra of the operators. $\hfill \square$

\section{The Calderon Problem}

While the edge calculus contains all “standard” elliptic BVPs for differential operators,
the result holds as long as there is no obstruction to ellipticity.
Obstruction to ellipticity is expressed as obstruction to the existence of Fredholm problems
for a given operator, or equivalently as obstruction to the existence of Fredholm operators with a given interior symbol.
It often happens that even though the interior symbol of an edge-degenerate
operator is elliptic, the edge symbol is not invertible and the operator fails to be
Fredholm. For the edge problem corresponding to the Calderon problem we will create a new edge symbol, by adding some “boundary and “coboundary” operators concentrated at the edge, which will offer to the edge symbol maximum domain of invertibility.

\subsection{Proof of Theorem \ref{Calderon}}
Finding out the cases for which the edge problem for the Calderon problem can be made Fredholm, indicates the cases for which the uniqueness argument of section \ref{proofell} applies. The following Lemma directly makes use of the analysis of \cite{SchEl}, Section 6.1.1 for Fredholmness of the Laplace operator, as it turns out the problem is essentially reduced to that of the Laplacian for $\sigma\ne 0$ and locally flat near the boundary manifold.

\begin{lem}
The operator Cald corresponding to the multidimensional Calderon problem on a manifold \( \mathcal{M} \)
$$
\text{Cald}\equiv\text{div}(\sigma\nabla):\mathcal{W}^{s,\gamma}(\mathcal{M})\rightarrow \mathcal{W}^{s-2,\gamma-2}(\mathcal{M}),
$$
where \( \sigma \) is the conductivity,  
can be made Fredholm through conversion to an edge problem via addition of boundary and coboundary operators, for 
$\gamma\ne \frac{1}{2},\frac{3}{2}$.
\end{lem}
\noindent {\it Proof}.

1. \textit{Operator expressed in local coordinates.} 
Consider local coordinates \( x=(r, x_1, x_2, \ldots, x_{n-1}) \) near the boundary, with \( r \) normal to the boundary and \( x_i \) tangential.
The Calderón operator in local coordinates is:
  \[
  \text{div}(\sigma(x) \nabla u) = \frac{1}{\sqrt{|g|}} \partial_i \left( \sqrt{|g|} g^{ij} \sigma(x) \frac{\partial u}{\partial x_j} \right)
  \]

Suppose the manifold is locally flat near the boundary, which means simplifying the metric tensor to approximately an identity matrix; hence, \( g_{ij} \approx \delta_{ij} \) and \( |g| \approx 1 \). Under this assumption, the Calderón operator can be simplified to:
  \[
  \text{div}(\sigma(x) \nabla u) \approx \frac{\partial}{\partial r} \left( \sigma(x) \frac{\partial u}{\partial r} \right) + \sum_{i=1}^{n-1} \frac{\partial}{\partial x_i} \left( \sigma(x) \frac{\partial u}{\partial x_i} \right)
  \]
where 
the divergence and gradient are reduced to their Euclidean forms.\\

2. \textit{Edge symbol of the operator}. 
The edge symbol for an operator is derived by applying a Fourier transform to the tangential derivatives of the operator in local coordinates.
As a pseudodifferential operator, the edge symbol is primarily concerned with capturing the leading-order behavior of the operator, particularly in a high-frequency domain; hence, as mixed derivatives do not contribute to the highest-order terms in this context, the focus is on pure higher-order derivatives.

Therefore, in the Fourier-transformed space, the edge symbol of the operator Cald becomes
  \[
  \sigma_{\wedge}(\text{Cald})(\xi) = \sigma(x)(\frac{\partial^2}{\partial r^2} - \|\xi\|^2),
  \]
where \(\|\xi\|^2 = \xi_1^2 + \xi_2^2 + \ldots + \xi_{n-1}^2\) is a multi-dimensional frequency variable, with \( \xi_i \) the Fourier dual variables corresponding to the tangential coordinates \( x_i \).

\textit{Hence, the foundational structure of the edge symbol corresponding to the Calderón problem retains a consistent format regardless of the specific form of conductivity \(\sigma(x)\), be it constant, linear, radial, polynomial, harmonic, exponential, etc.  }

3. \textit{Adjoint operator of the edge symbol.}  
For test functions \( f, g \in C_{0}^{\infty}(\mathbb{R}_{+} \times \mathbb{R}^{n-1}) \), we apply integration by parts, focusing on the radial derivative:
  \[
  \int_{\mathbb{R}_+ \times \mathbb{R}^{n-1}} \sigma(x) f \cdot \left(\frac{\partial^2 g}{\partial r^2} - \sum_{i=1}^{n-1}\xi_i^2 g\right) \, dr \, d\mathbf{x} = -\int_{\mathbb{R}_+ \times \mathbb{R}^{n-1}} \sigma(x) \frac{\partial f}{\partial r} \cdot \frac{\partial g}{\partial r} \, dr \, d\mathbf{x}
  \]

The adjoint calculation involves 'reversing' the differential operation in the radial direction while keeping the tangential frequency components:
  \[
  -\int_{\mathbb{R}_+ \times \mathbb{R}^{n-1}} \sigma(x) \frac{\partial f}{\partial r} \cdot \frac{\partial g}{\partial r} \, dr \, d\mathbf{x} = \int_{\mathbb{R}_+ \times \mathbb{R}^{n-1}} \sigma(x) \left(\frac{\partial^2 f}{\partial r^2} - \sum_{i=1}^{n-1}\xi_i^2 f\right) \cdot g \, dr \, d\mathbf{x}
  \]

This equation suggests that the adjoint family
$$
\sigma_{\wedge}(\text{Cald})(\xi)^{*}:\mathcal{K}^{2-s,2-\gamma}(\mathbb{R}_{+})\rightarrow \mathcal{K}^{-s,-\gamma}
(\mathbb{R}_{+}),
$$
on $C_{0}^{\infty}(\mathbb{R}_{+})$ retains the same structure as the original edge symbol $\sigma_{\wedge}(\text{Cald})(\xi)$.

4. \textit{Homogeneity of the edge symbol}. 
To study the invertibility of the edge symbol in the spaces
$$
\sigma_{\wedge}(\text{Cald})(\xi):\mathcal{K}^{s,\gamma}(\mathbb{R}_{+})\rightarrow \mathcal{K}^{s-2,\gamma-2}(\mathbb{R}_{+}),
$$
it suffices to study the invertibility of $\sigma_{\wedge}(\Delta)(\xi)$ at the normalized scale $|\xi|=1$, and then extend the result to all $\xi\neq 0$ across the entire domain by scaling, due to the principle of homogeneity of the edge symbol. 

5. \textit{Kernel and cokernel of the edge symbol.}
The kernel of the edge symbol is 
\[
\sigma_{\wedge}(\text{Cald})(\xi)u = \sigma(x)(\frac{\partial^2 u}{\partial r^2} - \sum_{i=1}^{n-1} \xi_i^2 u) = 0.
\] 
For the case of spatially varying conductivity \(\sigma(x)\ne 0\) we get the equation
\[
\frac{\partial^2 u}{\partial r^2} - \sum_{i=1}^{n-1} \xi_i^2 u = 0.
\]
whose general solution is:
\[ 
u(r, \xi_1, \xi_2, \ldots, \xi_{n-1}) = C_1 e^{r\sqrt{\sum_{i=1}^{n-1} \xi_i^2}} + C_2 e^{-r\sqrt{\sum_{i=1}^{n-1} \xi_i^2}}. 
\]

To be part of the kernel, the solution \(u\) must be integrable over \( \mathbb{R}_+ \). Given that the term \( e^{r\sqrt{\sum_{i=1}^{n-1} \xi_i^2}} \) grows exponentially as \( r \) increases (particularly for \( \sum_{i=1}^{n-1} \xi_i^2 \neq 0 \)), integrability requires that the coefficient \( C_1 \) associated with this term must be zero to prevent the function from becoming unbounded as \( r \rightarrow \infty \).
Therefore, the kernel for the Calderón problem's edge symbol 
has a one-dimensional kernel spanned by 
\begin{equation} \label{span}
u = e^{-|\xi|r}.
\end{equation}

Cokernel, i.e. the kernel of the adjoint edge symbol, $$
\sigma_{\wedge}(\Delta)(\xi)^{*}:\mathcal{K}^{2-s,2-\gamma}(\mathbb{R}_{+})\rightarrow \mathcal{K}^{-s,-\gamma}(\mathbb{R}_{+}),
$$
is also spanned by \(e^{-|\xi|r}\). 

The norm \(\mathcal{K}^{s,\gamma}(\mathbb{R}^n)\) of a function \( u \) is defined as
\[
\| u \|_{\mathcal{K}^{s,\gamma}} = \left( \sum_{|\alpha| \leq s} \int_{\mathbb{R}^n} r^{-2\gamma} |D^\alpha u|^2 \, dr \, d\xi \right)^{1/2},
\]
where \( D^\alpha \) represents derivatives up to order \( s \) with respect to both the radial variable \( r \) and the angular variables \(\xi\) within the multi-index \(\alpha\).
To avoid singularities at \( r = 0 \), the combined exponent on \( r \) from \( r^{-2\gamma} \) and the derivatives should be greater than \(-1\).
Hence, for the function \( u(r, \xi_1, \ldots, \xi_{n-1}) = e^{-r|\xi|} \), the integral's convergence relies on the term 
\[
\int_{\mathbb{R}^n} r^{-2\gamma} r^k |\xi|^k e^{-2r|\xi|} \, dr \, d\xi,
\]
where $k$ is the degree of \(r\) in \(D^\alpha u\). 
Therefore the condition \( k - 2\gamma > -1\) must hold for the smallest \( k \) (i.e., \( k = 0 \)), thus 
\[
\gamma < \frac{1}{2}.
\]

Hence, the presence of the weight factor $r^{-2\gamma}$ in the definition of the norm in
$\mathcal{K}^{s,\gamma}(\mathbb{R}_{+})$ results in the assertion that $e^{-r}\in \mathcal{K}^{s,\gamma}(\mathbb{R}_{+})$  
if and only if $\displaystyle \gamma<\frac{1}{2}$.

On the other hand, for the function \( u(r, \xi_1, \ldots, \xi_{n-1}) = e^{-r|\xi|} \), the convergence of $\| u \|_{\mathcal{K}^{2-s,2-\gamma}}$ relies on the term 
\[
\int_{\mathbb{R}^n} r^{-2(2-\gamma) + k} |\xi|^{2k} e^{-2r|\xi|} \, dr \, d\xi,
\]
where $k$ is the degree of \(r\) in \(D^\alpha u\). 
Therefore the condition $-2(2-\gamma) + k >-1$ must hold for the smallest \( k \) (i.e., \( k = 0 \)), thus 
\[
\gamma > \frac{3}{2}.
\]

Hence, the presence of the weight factor $r^{-2(2-\gamma)}$ in the definition of the norm in
$\mathcal{K}^{2-s,2-\gamma}(\mathbb{R}_{+})$ results in the assertion that $e^{-r|\xi|}\in \mathcal{K}^{2-s,2-\gamma}(\mathbb{R}_{+})$
if and only if $\displaystyle \gamma>\frac{3}{2}$.

As a result of the above, we conclude that:
\begin{itemize}
\item[\textbf{\textit{Case 1}}] If $\displaystyle \gamma<\frac{1}{2}$, then $\sigma_{\wedge}(\text{Cald})(\xi)$ has the one-dimensional kernel spanned by the
function (\ref{span}) and the trivial cokernel.

\item[\textbf{\textit{Case 2}}] If $\displaystyle \gamma>\frac{3}{2}$, then $\sigma_{\wedge}(\text{Cald})(\xi)$ has the trivial kernel and the one-dimensional
cokernel spanned by the function (\ref{span}).

\item[\textbf{\textit{Case 3}}] If $\displaystyle \frac{1}{2}<\gamma<\frac{3}{2}$, then $\sigma_{\wedge}(\text{Cald})(\xi)$ is invertible, satisfying both injectivity (trivial kernel) and surjectivity (trivial cokernel) conditions.

\item[\textbf{\textit{Case 4}}] If $\displaystyle \gamma=\frac{1}{2}$ or $\displaystyle \gamma=\frac{3}{2}$, then $\sigma_{\wedge}(\text{Cald})(\xi)$ is not
Fredholm. 
In particular, for \(\gamma = 1/2\) or (3/2), small alterations on functions belonging to \(\mathcal{K}^{s,\gamma}(\mathbb{R}^n)\) might lead to large deviations onto their images under the edge symbol, causing unexpected limit infinities outside the range; thus the requirement for a closed range is not satisfied.
\end{itemize}

6. We want to make our edge symbol invertible, i.e. make the equation 
\begin{equation} \label{SchEl_(6.2)}
\sigma_{\wedge}(\text{Cald})(\xi)v=F
\end{equation}
uniquely solvable by putting finitely many conditions on the solution or the right-hand side. This is obviously impossible for the non-closed range case 4 (because the limit of a converging sequence in the domain might map to an element outside the range, making it impossible to uniquely link elements in the target space back to the domain), and we will deal only with cases 1 and 2.

6a. \textit{Case 1}: Integral condition for Unique Solvability
for $\displaystyle \gamma<\frac{1}{2}$. 

The point is to restrict the solution space. For example, a global constraint over the whole solution space, determines the integral condition
\begin{equation} \label{SchEl_(6.3)}
B(\displaystyle \xi)v\equiv\int_{0}^{\infty}\phi(|\xi|r)v(r)dr=g\in \mathbb{C},
\end{equation}
where  $\phi(|\xi|r)\in C_{0}^{\infty}(\mathbb{R}{+})$ is not orthogonal to $e^{-r\sqrt{\sum{i=1}^{n-1} \xi_i^2}}$.

The non-orthogonality uniquely imposes a scalar constraint, ensuring that solutions \( v(r) \) are solely scalar multiples $c$ of the kernel function. 
Thus it guarantees the equation's unique solvability within the kernel by preventing any arbitrary function from trivially satisfying it, thereby directly linking each solution \(v(r)\) to a unique scalar \(c\) that corresponds with a given \(g\). Hence, the integral condition becomes a scalar equation for \( c \), and essentially effectively reduces the infinite-dimensional solution space to a singular dimension defined by a specific scalar \( c \in \mathbb{C} \). 

To show uniqueness, let \(v_1(r)\) and \(v_2(r)\) be two solutions of the integral condition for the same $g$. The difference \(v_d(r) = v_1(r) - v_2(r)\) is a solution of the homogeneous integral condition. As \(v_1(r)\) and \(v_2(r)\) are of the form of a scalar multiplied by the kernel function \(e^{-|\xi|r}\), then so is \(v_d(r)\). As \(\phi(r)\) is not orthogonal to the kernel function, then \(v_d(r)\) is identically zero, i.e. \(v_1(r) = v_2(r)\).

6a.i. \textit{Edge problem by quantization for Cald.}
The symbol of the Laplace-Beltrami operator is quantized by smoothing near \(\xi = 0\) (to ensure that the operator behaves well, even at low frequencies or near the "edges" of the domain) and 
replacing the multidimensional frequency vector \(\xi\) with the corresponding multidimensional differentiation operator \( -i\frac{\partial}{\partial x} \), obtaining the following edge problem for Cald operator for $\gamma<\frac{1}{2}$:
\begin{equation} \label{SchEl_(6.4)}
\left\{\begin{array}{l}
\text{Cald}\ u=f,\\
Bu\equiv\int_{0}^{\infty}\phi(r[-i\frac{\partial}{\partial x}])u(x)dr=g(x).
\end{array}\right.
\end{equation}

6a.ii. \textit{Estimates for Cald}.

In what spaces will this problem be Fredholm? The spaces in which $\text{Cald}$ acts are clear, and we only need to find the
natural space into which the edge boundary operator $B$ acts.
Namely, by Theorem 3.32 of \cite{SchEl} 
$$
B(\xi)\in S_{CV}^{0}(\bullet_{1},\ \bullet_{2}),
$$
where $S_{CV}^{0}$ is the space of symbols of compact fiber variation, and $\bullet_{1}$ is the space $\mathcal{K}^{s,\gamma}(\mathbb{R}_{+})$ equipped with the family of norms
$\Vert\cdot\Vert_{\xi}=\Vert|\xi|^s\kappa_{|\xi|}^{-1}\cdot\Vert_{\mathcal{K}^{s,\gamma}(\mathbb{R}_{\bullet})}
$
and $\bullet_{2}$ is the one-dimensional complex space $\mathbb{C}$ equipped with the family of norms
$\Vert\cdot\Vert_{\xi}=|\xi|^{s+1/2}|\cdot|.$

By \cite{SchEl}, Theorem 3.32,
the operator $B$ acts continuously in the spaces
$$
B\ :\ \mathcal{W}^{s,\gamma}(M)\rightarrow H^{s+1/2}(X),\ X=\partial M.
$$

6a.iii. \textit{Edge problem and Fredholmness for Cald}.
Hence, the operator A corresponding to our edge boundary value problem acts in the spaces
\begin{equation} \label{SchEl_(6.5)}
\mathrm{A}=\left(\begin{array}{l}
\Delta\\
B
\end{array}\right)\ :\mathcal{W}^{s,\gamma}(M)\rightarrow\ \begin{array}{c}
\mathcal{W}^{s-2,\gamma-2}(M) \\
\oplus \\
H^{s+1/2}(X)
\end{array}.
\end{equation}
This operator is Fredholm; this follows from the general finiteness theorem \cite{SchEl}, Theorem 6.19.

6b. \textit{Case 1}: Condition for Unique Solvability
for $\displaystyle \gamma>\frac{3}{2}$.

To make Eq. (\ref{SchEl_(6.2)}) uniquely solvable, one can equip it, say,
with a co-condition including a numerical unknown $\mu\in \mathbb{C}$:
\begin{equation} \label{SchEl_(6.6)}
\sigma_{\wedge}(\Delta)(\xi)v+\mu\phi(|\xi|r)=F,
\end{equation}
where $\phi$ is the same function as above.

6b.i. 
Quantizing, we obtain the problem
\begin{equation} \label{SchEl_(6.7)}
\Delta u+Cw=f,
\end{equation}
where the operator $C$ is given by the formula
$$
Cw=\phi\left(r\left[-i\frac{\partial}{\partial x}\right]\right)w.
$$

6b.ii.  
A similar argument shows that the operator A corresponding to this edge coboundary value problem acts in the spaces
\begin{equation} \label{SchEl_(6.8)}
\mathrm{A}= (\Delta\ C) : \begin{array}{c}
\mathcal{W}^{s,\gamma}(M) \\
\oplus \\
H^{s-5/2}(X)
\end{array}\rightarrow \mathcal{W}^{s-2,\gamma-2}(M).
\end{equation}
$\hfill \square$

\section*{Appendix:  Overview of Schulze's  Pseudodifferential edge algebra}

\textbf{Construction of admissible algebras in manifolds with singularities} \cite{SchEl}.
A singular manifold $(M^{\mathrm{o}},\ \mathcal{D})$ is defined by its differential operators' behavior, where $M^{\mathrm{o}}$ is a smooth manifold and $\mathcal{D}$ an algebra of differential operators. These operators are standard in $M^{\mathrm{o}}$ but adhere to specific constraints near singular points, ensuring the algebra $\mathcal{D}$ within any compact subset $U\Subset M^{\mathrm{o}}$ aligns with all differential operators having smooth coefficients on $M^{\mathrm{o}}$ \cite{Qm, SSSa}. This framework facilitates handling differential equations on manifolds with singularities, focusing on operator behavior rather than the manifold's embedding or metric properties. The distinction among singular manifolds sharing the same $M^{\mathrm{o}}$ lies in the differential operators' limits at the manifold's infinite regions, as outlined in \cite{SSSa}.

\textit{Algebras of differential operators on manifolds with cone and edge singularities.}\cite{SchEl}, primarily Section 1.1.2.
The algebra $\mathcal{D}$ for differential operators on manifolds with singularities or edges is generated by $C^{\infty}(M)$ and a space $V$ of vector fields, where $V$ and $F$, a function space, have inclusions with respect to smooth functions and vector fields on $M^{\mathrm{o}}$. By embedding $M^{\mathrm{o}}$ into a compact manifold $M$ and extension of a non-degenerate Riemannian metric $d\rho^2$ from $M^{\mathrm{o}}$ to $M$, $V$ is then characterized via its dual space $V'$ based on an $F$-valued inner product, making $V=V'$. This construction uniquely defines $\mathcal{D}$, which then ties manifolds to two categories; with conical singularities or edges, depending on the operator type.

Symbols for cone-degenerate differential operators. 
\cite{SchEl}, Section 1.2.6. 
The principal symbol, analogous to that in smooth manifolds, emerges within the structure of the differential operator algebra $\mathcal{D}$, organized by order into a hierarchy $\mathcal{D}_k$. This organization leads to a graded algebra $\mathrm{gr}\ \mathcal{D} = \bigoplus_{j=0}^{\infty} \mathcal{D}_{j}/\mathcal{D}_{j-1}$, with $\mathcal{D}_{-1} = \{0\}$. For an operator $D$ of order $k$, its principal symbol $\Sigma(D) = \Sigma_{k}(D)$ is identified through projection to $\mathrm{g}\mathrm{r}_{k}\mathcal{D}$. In smooth contexts, this symbol functions over the cotangent bundle $T^{*}M$, guiding from space $M$ to $T^{*}M$. The transition, or microlocalization, especially for manifolds with conical singularities, lacks a clear phase space definition, prompting a two-stage process to navigate potential complexities.

Stage 1. Localization. \cite{SchEl}. 
Functions in $C^{\infty}(\mathcal{M})$, constant on the boundary, are a subalgebra in the differential operator algebra $\mathcal{D}$, exhibiting a commutation relation indicative of a $C^{\infty}(\mathcal{M})$-module structure in the associated graded algebra $\mathrm{gr}\mathcal{D}$. Symbol spaces $\Sigma_{kx}$ localize the symbol order $k$ to a point $x$ on $\mathcal{M}$. The set of local representatives of an operator $D$ in $\mathcal{D}_{k}$, denoted $\sigma_{x}(D)$, uniquely determines $D$'s principal symbol. These local representatives can be conveniently described via a scaling procedure related to the metric's induced local scaling transformations $g_{\lambda}$ around point $x$. This leads to defining the scaled operator $D_{x}$ by the limit of scaled operations on $D$, with the limit interpreted in terms of pointwise convergence. This limit exists for any $D$ in $\mathcal{D}_{m}$ and establishes a bijective, multiplicative correspondence between $\sigma_{x}(D)$ and $D_{x}$. This localizes each operator $D$ to a family of operators $D_{x}$, representing $D$'s class in $\mathrm{gr}_{m}\mathcal{D}$ on the space $K_{x}$, adjusted for interior points and the conical point specifically.

Stage 2. Microlocalization. \cite{SchEl}.
For each interior point $x$, the constant-coefficients operator $D_x$ on $T_{x}M$ leads to the polynomial $P_{m}(0,\ \xi)$ via Fourier transform, representing $D$'s principal symbol in $T^{*}M$'s fiber over $x$. This symbol, $\sigma(D)$, extends to a smooth function across $T^{*}\mathcal{M}$, defining the {\it interior symbol}. Conversely, at the conical point, $D_{x}$ lacks translation invariance, precluding microlocalization, and is directly termed the {\it cone symbol}, $\sigma_{c}(D)$, preserving the distinct treatments for interior and conical points in symbol analysis.

Symbols for edge-degenerate differential operators.\cite{SchEl}. 
The algebra $\mathcal{D}$ of edge-degenerate differential operators is filtered by operator order, defining symbols within  $\mathrm{g}\mathrm{r}\mathcal{D}=\bigoplus_{m=0}^{\infty}\mathcal{D}_{m}/\mathcal{D}_{m-1}$. The principal symbol of an operator $P\in \mathcal{D}_{m}$ is its image in $\mathcal{D}_{m}/\mathcal{D}_{m-1}$. This principal symbol, detailed for edge-degenerate operators, combines the {\it interior symbol} and the {\it edge symbol}. The interior symbol, $\sigma(D)$, derived from the classical symbol $\sigma_{clas}(D)$, extends smoothly to $T^{*}\mathcal{M}$. The edge symbol, $\sigma_{\wedge}(D)(x,\ \xi)$, represented in local coordinates, forms a well-defined operator family on $K_{\Omega_{x}}$, parameterized by $T^{*}X$ points. The compatibility condition ensures the interior and edge symbols of $D\in \mathcal{D}_{m}$ satisfy $\sigma(\sigma_{\wedge}(D))=\sigma_{\partial}(D)$, integral to symbol definitions for edge-degenerate differential operators on manifolds with edges. This setup indicates that $M$ acts as the stretched manifold of $\mathcal{M}$.

\textit{Algebras of pseudodifferential operators on manifolds with cone and edge singularities.}
The task involves extending the algebra $\mathcal{D}_{k}(M)$ to $PS\mathcal{D}_{k}(M)$, including pseudodifferential operators, with emphasis on handling singular points on $\mathcal{M}$. The approach is local, constructing these operators in coordinate charts and integrating them across the manifold using partitions of unity, similar to conventional pseudodifferential operator theory \cite{[H4]}. Special focus is on constructing admissible  pseudodifferential operators near singularities, extending symbol classes beyond polynomials via noncommutative analysis \cite{Qm}.

STEP 1.\textit{Function spaces.}
The algebra for pseudodifferential operators is defined within a scale of Hilbert spaces, essential for ensuring the operators act in function spaces suitable for singular elliptic differential operators \cite{Qm}. These spaces must align with solutions of differential equations, with the choice proven critical for theoretical integrity \cite{SSSa}, \cite{Qm}. Specifically, within 
$\mathcal{M}$'s interior, these function spaces match traditional Sobolev spaces, while near singularities, the spaces adapt to the singularity type, differentiated into cone-degenerate and edge-degenerate categories \cite{SchEl}.

Function spaces for cone-degenerate operators. \cite{SchEl}, Section 2.1.1.
Weighted Sobolev spaces $H^{s,\gamma}(\mathcal{M})$ are tailored for cone-degenerate elliptic operators, incorporating solutions of homogeneous equations and exhibiting a norm invariance property near singularities \cite{SchEl}. These spaces match ordinary Sobolev spaces in $\mathcal{M}$'s interior and adapt to singularity types near conical points. The norm calculation for $H^{s,\gamma}(\mathcal{M})$ relies on functions near singularities, transitioning to the model cone $K=K_{\Omega}$ for a global group action. The weighted Sobolev space norm combines standard Sobolev norms and localized norms near singularities, forming a scale of Hilbert spaces invariant under specific transformations \cite{SchEl}. Cone-degenerate differential operators maintain continuity within these spaces, ensuring operator applicability across different Sobolev scales \cite{SchEl}.

Function spaces for edge-degenerate operators. \cite{SchEl}, Section 2.1.1. 
In studying edge-degenerate operators, we define function spaces that accommodate the operators' behavior near and away from the manifold's edge, ensuring coherence across these regions. For the cone $K_{\Omega}$, we integrate usual Sobolev spaces at infinity with weighted Sobolev spaces near the vertex, employing a partition of unity for seamless transitioning. This approach creates the space $\mathcal{K}^{s,\gamma}(K_{\Omega})$, blending local and global function norms (\cite{PsMSi}).
Abstract wedge spaces $\mathcal{W}^{s}(\mathbb{R}^{n},\ H)$ extend this concept to manifold edges, adapting through Fourier transform methods and scaling by a strongly continuous operator group $\kappa_{\lambda}$, leading to the specific space $\mathcal{W}^{s,\gamma}(W)$ for an infinite wedge. This construction is then applied to the entire manifold $\mathcal{M}$, forming $\mathcal{W}^{s,\gamma}(\mathcal{M})$ which integrates edge behavior with the manifold's bulk properties.
An edge-degenerate differential operator $D$ operates continuously within these constructed spaces, confirming the adequacy of the spaces $\mathcal{W}^{s,\gamma}(\mathcal{M})$ for hosting solutions to edge-degenerate differential equations \cite{SchEl}.

STEP 2. \textit{Symbols}. \textit{Describe the admissible class of symbols.} \cite{Qm}.
For applications in noncommutative analysis, the symbol class $S^{\infty}(\mathbb{R}^{n})$ is optimal. Yet, in the context of pseudodifferential operators, symbols typically have uniform dependence on spatial variables without infinite growth. Symbols that do exhibit growth are not bounded within the relevant spaces, nor are they part of the algebra, a fact that applies to manifolds with singularities too. To facilitate analysis of such operator algebras and symbol classes, the approach to functions of noncommuting operators is expanded to include operator-valued symbols.

Symbols for cone-degenerate pseudodifferential operators.\cite{SchEl}, Section 3.4.1.
On a manifold $\mathcal{M}$ with conical singularities, the principal symbol of a cone-degenerate pseudodifferential operator $D$ combines the interior symbol $\sigma(D)$ and the cone symbol $\sigma_{c}(D)$. The interior symbol, a function on the stretched cotangent bundle $T_{0}^{*}\mathcal{M}$, is homogeneous of order $m$ in the fibers. The cone symbol, defined through the concept of conormal symbols, is an mth-order pseudodifferential operator with parameter $p$ on the base $\Omega$ of the cone, continuous in weighted Sobolev spaces. Compatibility between interior and cone symbols is required for them to form a principal symbol, indicated when the principal symbol of the cone symbol matches the restriction of the interior symbol to $\partial T_{0}^{*}\mathcal{M}$. This framework accommodates the analysis of pseudodifferential operators on manifolds with singularities, encapsulating both local and global operator characteristics \cite{AV}, \cite{SchEl}.

Symbols for cone-degenerate pseudodifferential operators.\cite{SchEl}, Section 3.4.1.
On a manifold $\mathcal{M}$ with conical singularities, the principal symbol of a cone-degenerate pseudodifferential operator $D$ combines the interior symbol $\sigma(D)$ and the cone symbol $\sigma_{c}(D)$. The interior symbol $\sigma=\sigma(y,\ \eta)$ is a function on the stretched cotangent bundle $T_{0}^{*}\mathcal{M}$, homogeneous of order $m$ in the fibers. Cone symbols are defined via conormal symbols, which are mth-order pseudodifferential operators on the manifold $\Omega$ with a specific parameter. These symbols, when associated with a conormal symbol, form the cone symbol $\sigma_{c}$, operating in weighted Sobolev spaces on the infinite cone $K_{\Omega}$. Compatibility between $\sigma$ and $\sigma_{c}$ ensures that they form principal symbols of the operator $D$, characterized by a shared boundary symbol $\sigma_{\partial}$ \cite{AV}.

Symbols for edge-degenerate pseudodifferential operators. \cite{SchEl}, Section 3.4.2. 
For compact manifolds $\mathcal{M}$ with an edge $X$ and a cone base $\Omega$, the principal symbols of edge-degenerate operators consist of interior symbols $\sigma$ and edge symbols $\sigma_{\wedge}$. Interior symbols are homogeneous functions on the stretched cotangent bundle, while edge symbols represent families of differential or pseudodifferential operators on $K_{\Omega}$, parameterized by $T_{0}^{*}X$ and acting in $\mathcal{K}^{s,\gamma}(K_{\Omega})$. Edge symbols incorporate conditions for twisted homogeneity, continuity, almost compact fiber variation, and are defined in both local neighborhoods of the cone vertex and its exterior. The compatibility between interior symbols and cone symbols ensures the principal symbol's coherence. Edge symbols, represented by the formula $\sigma_{c}=r^{-m}\sigma_{c}\left(ir\frac{\partial}{\partial r}\right)$, are continuous in weighted Sobolev spaces, signifying the structured approach in defining operator families for edge-degenerate pseudodifferential operators on manifolds with singular edges.

STEP 3. \textit{Quantization} 

Quantization on manifolds with cones. \cite{SchEl}, Section 3.4.1. 
In order to thus associate a pseudodifferential operator in weighted Sobolev spaces on on $\mathcal{M}$ to each principal symbol, we introduce negligible operators, which are compact in $D:H^{s,\gamma }( \mathcal{M} )\rightarrow H^{s-m,\gamma -m}( \mathcal{M} )$, and continuous for $D:H^{s,\gamma }( \mathcal{M} )\rightarrow H^{s-m+1,\gamma -m}( \mathcal{M} )$. We then construct a pseudodifferential operator $D$ of order $m$ and weight $\gamma$ with principal symbol $\bullet$ to be continuous in these spaces, its definition modulo negligible operators. These operators, denoted by $\Psi_{\gamma}^{m}(M)$, have well-defined principal symbols, and the principal symbol $\bullet(D)$, uniquely characterizes each $D\in\Psi_{\gamma}^{m}(M)$. The compactness of it within $\Psi_{\gamma}^{m}(M)$ signals its location in the space of negligible operators $I_{\gamma}^{m}$.

Quantization on manifolds with edges. \cite{SchEl}, Section 3.4.3.
Operators of order $m$ on the compact manifold $\mathcal{M}$ with edges are defined as those extendable to continuous from $\mathcal{W}^{s,\gamma}(\mathcal{M})$ to $\mathcal{W}^{s-m,\gamma-m}(\mathcal{M})$ for every $s\in \mathbb{R}$. These are called {\it negligible operators} when compactly-mapped between these spaces and are included in $J\mathrm{O}\mathrm{p}{\gamma}^{m}(\mathcal{M})$. The pseudodifferential operators of order $m$ and weight $\gamma$ have unique principal symbols, which can be written down as composed of interior symbol $\sigma\in C^{\infty}(T{0}^{*}\mathcal{M})$ and edge symbol $\sigma_{\wedge}$. These symbols need to satisfy compatibility condition. Then they can be together quantized, and the obtained pair $(\sigma,\ \sigma_{\wedge})$ is unique for each operator $P$. This ensures that the operator is precisely defined from its principal symbol within the framework of weighted Sobolev spaces on a compact manifold with edges. The resulting pseudodifferential operators’ algebra is denoted by $\mathrm{P}\mathrm{S}\mathrm{D}_{\gamma}^{m}(\mathcal{M})$.

STEP 4. \textit{Algebra formation}. {\it We study whether the symbols of the aforementioned pseudodifferential operators form an algebra.} 
In \cite{SchEl}, the algebra of principal symbols for zero-order pseudodifferential operators ($\psi\mathrm{D}\mathrm{O}$s) (the general case follows
by order reduction) is defined within a Hilbert space $\mathcal{H}$, emphasizing the role of a symbol mapping $\sigma$ that distinguishes operators by their association to the ideal $\mathcal{K}$ of compact operators. This setup yields a monomorphism $\tilde{\sigma}$ of the Calkin algebra into a unital topological algebra $S$, alongside a linear mapping $Q$ serving as the right inverse of $\sigma$, facilitating quantization. Essentially, $S$ represents the algebra of (principal) symbols for $\mathcal{\tilde{A}}$, making it isomorphic to the corresponding Calkin algebra, which is a $C^*$ algebra.

\textbf{Edge problems}. \cite{SchEl}, Chapter 6.
In the theory of operators on manifolds with edges, boundary value problems are a specific instance where the edge acts as the boundary and the cone simplifies to $\mathbb{R}_+$, representing the inner normal \cite{PsMSi}. Addressing elliptic problems on these manifolds involves modifying the problem to ensure the edge symbol becomes invertible, often through the inclusion of boundary and coboundary operators in a matrix operator setup. This adaptation, rooted in the physical context of boundary conditions in Fredholm problems, may not always be feasible due to potential topological constraints akin to the Atiyah-Bott obstruction in boundary value problem theory.

\textit{Edge boundary and coboundary operators}. \cite{SchEl}, Section 6.1.2. 
In addressing Fredholm problems for operators on manifolds with edges, matrix operators incorporating edge boundary and coboundary conditions are essential. These operators, organized in a $2\times 2$ matrix format, integrate pseudodifferential operators with specific edge-focused components. This approach, critical for ellipticity and parabolicity studies, involves operators that are defined modulo negligible ones, aiming for operators that are compact and continuous across specified weighted Sobolev spaces \cite{SchEl}. The framework extends to defining both edge boundary and coboundary symbols as pseudodifferential operators on $X$, emphasizing their compact and continuous nature in adapted Sobolev spaces, ensuring the operators are finely tuned to the edge's properties. This method underscores a comprehensive strategy for handling complex operator conditions on manifolds with edges.

\textit{The calculus of edge morphisms.} \cite{SchEl}, Section 6.1.3.
Defining edge boundary and coboundary operators within a framework for operators on manifolds with edges, \cite{SchEl} introduces a comprehensive classification of operators as morphisms of order $m$ and weight $\gamma$. These morphisms, formed in a matrix operator setup combining pseudodifferential, boundary, and coboundary components, allow for the representation of complex boundary conditions integral to Fredholm problems. Key to this formalism is the concept of principal symbols for these morphisms, encapsulating both the interior symbol of the operator $A$ and an operator family forming the edge symbol, characterized by continuity and compactness properties in weighted Sobolev spaces. This structure underpins the analytical approach to elliptic problems on manifolds with edges, emphasizing compatibility conditions essential for morphism definition and the subsequent mathematical treatment of boundary phenomena.

\textbf{Inverse elliptic BVPs via BC method}

\textit{Ellipticity and finiteness theorems}

Ellipticity and finiteness theorems for degenerate $\psi \mathrm{D}\mathrm{O}$s. \cite{SchEl}, Section 3.5.1. and Section 3.5.4. 
In the degenerate $\psi \mathrm{D}\mathrm{O}$ algebra $\mathcal{A}$ constructed by \cite{SchEl}, an operator $A$ is termed elliptic if its symbol $\sigma(A)$ is invertible in the symbol algebra $S$. The main result of elliptic theory, the finiteness theorem, asserts that an operator is Fredholm if its symbol is invertible. This theorem applies to cone-degenerate and edge-degenerate $\psi \mathrm{D}\mathrm{O}$s on manifolds with conical singularities and edges, respectively, defining ellipticity through the invertibility of operators' symbols and stating that elliptic operators are Fredholm, with their kernel, cokernel, and index remaining constant across different smoothness levels.

Ellipticity and finiteness theorems for edge problems. \cite{SchEl}, Section 6.1.4. 
For a morphism $A$ within the degenerate $\psi \mathrm{D}\mathrm{O}$ algebra $\mathcal{A}$, \cite{SchEl} defines it as elliptic if its interior symbol is invertible across the entire cotangent bundle $T_{0}^{*}\mathcal{M}$ and its edge symbol is invertible on $T_{0}^{*}X$. Such elliptic morphisms are Fredholm across all considered spaces, maintaining constant kernel, cokernel, and index regardless of smoothness level $s$. Additionally, an elliptic morphism is uniquely identified by its principal symbol modulo negligible operators, underlining the pivotal role of the symbol in determining ellipticity and Fredholm properties.

The obstruction to ellipticity. \cite{SchEl}, Section 6.2. 
For an operator $A$ on a manifold with an edge, constructing an elliptic edge problem requires adding boundary operators $B$, coboundary operators $C$, and a pseudodifferential operator $D$ on the edge. The solvability of such a problem hinges on the interior ellipticity of $A$ and a Fredholm condition on its edge symbol, as detailed in \cite{SchEl}. The crux is that for $A$ to be part of an elliptic edge problem, its edge symbol's index must vanish in a specific $K$-theory element, influenced solely by $A$'s interior symbol.

\textbf{Inverse parabolic BVPs via BC method}. \cite{p}. 
In the complex upper half-plane $\mathbb{H}$, we explore symbol spaces $S^{\mu;\ell}(\mathbb{R}^{n}\times \mathbb{R}^{q};E,\tilde{E})$ composed of smooth, bounded functions from $\mathbb{R}^{n}\times \mathbb{R}^{q}$ into $\mathcal{L}(E,\tilde{E})$, where $E$ and $\tilde{E}$ are Hilbert spaces. These spaces, equipped with a topology defined by a seminorm system, include symbols that are (anisotropic) homogeneous of degree $\mu$. We also define the space of symbols of infinite order, $S^{-\infty}(\mathbb{R}^{n}\times \mathbb{R}^{q};E,\tilde{E})$, and extend these definitions to symbols dependent on spatial variables. Furthermore, classical symbols $S_{\mathrm{cl}}^{\mu;\ell}(\mathbb{R}^{n}\times \mathbb{R}^{q};E,\tilde{E})$ comprise a subset characterized by an asymptotic sum of homogeneous components. Similarly, for the Volterra property in symbols, we define $S_{V}^{\mu;\ell}(\mathbb{R}^{n}\times \mathbb{H};E,\tilde{E})$ emphasizing operators analytically extendable into the complex plane's interior, alongside corresponding definitions for pseudodifferential operators. Parabolic symbols and operators within this framework are distinguished by their parameter-dependent ellipticity. This detailed symbolic structure underpins the analysis of operators on manifolds, particularly emphasizing the treatment of elliptic and parabolic types in various contexts. Essential to this theory is the notion of parabolicity for operators in manifolds with edges, as encapsulated in \cite{SchInd}, where both the interior symbol's parabolic nature and the edge symbol's invertibility criterion play crucial roles.


\begin{thebibliography}{99}

\bibitem{LSch}
Liu, Xiaochun and Schulze, Bert-Wolfgang, \emph{Boundary value problems with global projection conditions}, Springer International Publishing, 2018.

\bibitem{SchInd}
Rempel, Stephan and Schulze, B.-W., ``Index theory of elliptic boundary problems,'' \emph{Mathematische Monographien} \textbf{55}, 1982.

\bibitem{SchBPS}
Schulze, B.-W., \emph{Boundary value problems and singular pseudo-differential operators}, Interscience Series of Texts, Monographs, and Tracks, 1998.

\bibitem{sing}
Nazaikinskii, V. E., Savin, A. Y., Schulze, B. W. and Sternin, B. Y., \emph{Elliptic theory on singular manifolds}, Chapman and Hall/CRC, 2005.

\bibitem{mix}
Harutyunyan, Gohar and Schulze, Bert-Wolfgang, \emph{Elliptic Mixed, Transmission and Singular Crack Problems}, European Mathematical Society, vol.~4, 2007.

\bibitem{[10]}
de Monvel, L. B., ``Boundary problems for pseudo-differential operators,'' \emph{Acta Mathematica} \textbf{126}(1): 11--51, 1971.

\bibitem{[4]}
Atiyah, M. F., Patodi, V. K. and Singer, I. M., ``Spectral asymmetry and Riemannian geometry. I.,'' \emph{Mathematical Proceedings of the Cambridge Philosophical Society} \textbf{77}(1): 43--69, January 1975, Cambridge University Press.

\bibitem{[5]}
Atiyah, M. F., Patodi, V. K. and Singer, I. M., ``Spectral asymmetry and Riemannian geometry. II.,'' \emph{Mathematical Proceedings of the Cambridge Philosophical Society} \textbf{78}(3): 405--432, November 1975, Cambridge University Press.

\bibitem{[6]}
Atiyah, M. F., Patodi, V. K. and Singer, I. M., ``Spectral asymmetry and Riemannian geometry. III.,'' \emph{Mathematical Proceedings of the Cambridge Philosophical Society} \textbf{79}(1): 71--99, January 1976, Cambridge University Press.

\bibitem{[Kr51a]}
Krein, M. G., ``Determination of the density of an inhomogeneous string from its spectrum,'' \emph{Doklady Akad. Nauk. SSSR} \textbf{76}: 345--348, 1951 (in Russian).

\bibitem{[Kr51b]}
Krein, M. G., ``On inverse problems for an inhomogeneous string,'' \emph{Doklady Akad. Nauk. SSSR} \textbf{82}: 669--672, 1951 (in Russian).

\bibitem{[Bla71a]}
Blagovestcenskii, A., ``The local method of solution of the non-stationary inverse scattering problem for an inhomogeneous string,'' \emph{Trudy Mat. Inst. Steklova} \textbf{115}: 28--38, 1971 (in Russian).

\bibitem{[Be87]}
Belishev, M. I., ``An approach to multidimensional inverse problems for the wave equation,'' \emph{Dolkl. Akad. Nauk SSSR} \textbf{297}: 524--527, 1987 (Engl. Transl.: Soviet Math. Dokl. 36 (1988), 481--484).

\bibitem{belku}
Belishev, M. I. and Kurylev, Ya V., ``Nonsteady inverse problem for the multidimensional wave equation "in the large",'' \emph{Journal of Soviet Mathematics} \textbf{50}(6): 1944--1951, 1990.

\bibitem{[KKL01]}
Kachalov, Alexander, Kurylev, Yaroslav and Lassas, Matti, \emph{Inverse boundary spectral problems}, Chapman and Hall/CRC, 2001.

\bibitem{[Bla71b]}
Blagovestcenskii, A., ``The nonselfadjoint inverse matrix boundary problem for a hyperbolic differential equation,'' in \emph{Problems of mathematical physics, 5, Spectral Theory}, Izdat. Leningrad Univ., Leningrad, 38--62, 1971 (in Russian).

\bibitem{[BeBla92]}
Belishev, M. I. and Blagovestchenski, A. S., ``Multidimensional analogs of the Gel’fand–Levitan–Krein equations in inverse problems for the wave equation,'' in \emph{Ill-Posed Problems of Mathematical Physics and Analysis}, Nobosibirsk: Nauka, 50--63, 1992 (in Russian).

\bibitem{IsKu}
Isozaki, Hiroshi and Kurylev, Yaroslav, \emph{Introduction to Spectral Theory and Inverse Problem on Asymptotically Hyperbolic Manifolds}, arXiv preprint \url{arXiv:1102.5382}, 2011.

\bibitem{bel17_[22]}
Belishev, Mikhail I., ``The Calderon problem for two-dimensional manifolds by the BC-method,'' \emph{SIAM J. Math. Anal.} \textbf{35}(1): 172--182, 2003.

\bibitem{bel17_[62]}
Belishev, Mikhail I. and Wada, Naoki, ``A C*-algebra associated with dynamics on a graph of strings,'' \emph{Journal of the Mathematical Society of Japan} \textbf{67}(3): 1239--1274, 2015.

\bibitem{bel17_[41]}
Belishev, M. I. and Demchenko, M. N., ``Elements of noncommutative geometry in inverse problems on manifolds,'' \emph{Journal of Geometry and Physics} \textbf{78}: 29--47, 2014.

\bibitem{bel17_[42]}
Belishev, M., Demchenko, M. and Popov, A., ``Noncommutative geometry and the tomography of manifolds,'' \emph{Transactions of the Moscow Mathematical Society} \textbf{75}: 133--149, 2014.

\bibitem{bel17_[37]}
Belishev, Mikhail Igorevich, ``On algebras of three-dimensional quaternionic harmonic fields,'' arXiv preprint \url{arXiv:1611.08523}, 2016.

\bibitem{bel17_[34]}
Belishev, M. I., ``A unitary invariant of a semi-bounded operator in reconstruction of manifolds,'' \emph{Journal of Operator Theory}: 299--326, 2013.

\bibitem{bel17_[56]}
Belishev, M. and Simonov, S., ``Wave model of the Sturm–Liouville operator on the half-line,'' \emph{St. Petersburg Mathematical Journal} \textbf{29}(2): 227--248, 2018.

\bibitem{bel17}
Belishev, Mikhail I., ``Boundary control and tomography of Riemannian manifolds (the BC-method),'' \emph{Russian Mathematical Surveys} \textbf{72}(4): 581, 2017.

\bibitem{Qm}
Nazaikinskii, V. E., Schulze, B.-W. and Sternin, B. Y., \emph{Quantization Methods in the Theory of Differential Equations}, CRC Press, 2002.

\bibitem{hee0}
Heunen, Chris and others, ``A topos for algebraic quantum theory,'' \emph{Comm. Math. Phys.} \textbf{291}: 63--110, 2009.

\bibitem{hee}
Heunen, Chris and others, ``The Gelfand spectrum of a noncommutative C*-algebra: a topos-theoretic approach,'' \emph{Journal of the Australian Mathematical Society} \textbf{90}(1): 39--52, 2011.

\bibitem{DosS}
Dos Santos Ferreira, David and others, ``The linearized Calderón problem in transversally anisotropic geometries,'' \emph{International Mathematics Research Notices} \textbf{2020}(22): 8729--8765, 2020.

\bibitem{bal}
Bal, Guillaume, ``Mathematical Tools for Imaging and Inverse Problems,'' 2024. Available at \url{https://www.stat.uchicago.edu/~guillaumebal/PAPERS/Mathematical_Tools.pdf}.

\bibitem{SchEl}
Nazaikinskii, V. E., Savin, A. Y., Schulze, B. W. and Sternin, B. Y., ``Elliptic theory on singular manifolds,'' Chapman and Hall/CRC, 2005.

\bibitem{SSSa}
Schulze, B.-W., Sternin, B. and Shatalov, V., \emph{Differential equations on singular manifolds: Semiclassical theory and operator algebras}, vol.~15, Wiley-VCH, 1998.

\bibitem{[H4]}
Hörmander, L., \emph{The analysis of linear partial differential operators}, Vols. 1--4, Springer-Verlag, New York, 1983/85.

\bibitem{PsMSi}
Schulze, B.-W., \emph{Pseudo-differential operators on manifolds with singularities}, Elsevier, 1991.

\bibitem{AV}
Agranovich, Mikhail Semenovich and Vishik, Marko Iosifovich, ``Elliptic problems with a parameter and parabolic problems of general type,'' \emph{Uspekhi Matematicheskikh Nauk} \textbf{19}(3): 53--161, 1964.

\bibitem{p}
Krainer, Thomas and Schulze, Bert-Wolfgang, ``On the inverse of parabolic systems of partial differential equations of general form in an infinite space-time cylinder,'' in \emph{Parabolicity, Volterra Calculus, and Conical Singularities}, Birkhäuser, 93--278, 2002.

\end{thebibliography}
\end{document}